%% file: main.tex
\documentclass[ijoc,nonblindrev]{informs3_2} 
\OneAndAHalfSpacedXII % current default line spacing
\usepackage{natbib}
\bibpunct[, ]{(}{)}{,}{a}{}{,}%
\def\bibfont{\small}%
\def\BIBand{and}%

\TheoremsNumberedThrough     % Preferred (Theorem 1, Lemma 1, Theorem 2)
\EquationsNumberedThrough    % Default: (1), (2), ...
\usepackage{subfig, float}
\usepackage{multirow, makecell}
\usepackage[dvipsnames]{xcolor}
\usepackage[prependcaption,textsize=small]{todonotes}

\newcommand{\Pa}{\widetilde{P}}
\newcommand{\U}{\mathcal{U}}
\newcommand{\Q}{\mathcal Q}
\newcommand{\C}{\mathcal{C}}
\newcommand{\D}{\mathcal{D}}
\newcommand{\G}{\mathcal{G}}
\newcommand{\N}{\mathcal{N}}
\newcommand{\E}{\mathcal{E}}
\newcommand{\Gc}{\widetilde{G}}
\newcommand{\Ec}{\widetilde{E}}
\newcommand{\Nc}{\widetilde{N}}
\newcommand{\Sn}{\mathcal{S}}
\newcommand{\Ph}{\mathcal P}
\newcommand{\Real}{\mathbb{R}}

\begin{document}
	%%%%%%%%%%%%%%%%
	
	%\setcounter{page}{1}
	%\DOI{10.1287/xxxx.0000.0000}%
	
	\RUNAUTHOR{Byeon, Van Hentenryck, Bent, and Nagarajan}
	
	\RUNTITLE{Communication-Constrained Expansion Planning}
	
	% \TITLE{Bundling Information Goods of Decreasing Value}
	%Enter the full title:
	\TITLE{Communication-Constrained Expansion Planning \\ for Resilient Distribution Systems}
	
	\ARTICLEAUTHORS{%
		\AUTHOR{Geunyeong Byeon and Pascal Van Hentenryck}
		\AFF{Industrial and Operations Engineering, University of Michigan \EMAIL{} \URL{}}
		\AUTHOR{Russell Bent and Harsha Nagarajan}
		\AFF{Los Alamos National Laboratory \EMAIL{} \URL{}}
	} % end of the block
	
	\ABSTRACT{%
          Distributed generation and remotely controlled switches have
          emerged as important technologies to improve the resiliency
          of distribution grids against extreme weather-related
          disturbances. Therefore it becomes important to study how
          best to place them on the grid in order to meet a
          resiliency criteria, while minimizing costs and capturing
          their dependencies on the associated communication systems
          that sustains their distributed operations. This paper
          introduces the Optimal Resilient Design Problem for
          Distribution and Communication Systems (ORDPDC) to address
          this need. The ORDPDC is formulated as a two-stage
          stochastic mixed-integer program that captures the physical
          laws of distribution systems, the communication connectivity
          of the smart grid components, and a set of scenarios which
          specifies which components are affected by potential
          disasters. The paper proposes an exact branch-and-price
          algorithm for the ORDPDC which features a strong lower bound
          and a variety of acceleration schemes to address degeneracy.
          The ORDPDC model and branch-and-price algorithm were
          evaluated on a variety of test cases with varying disaster
          intensities and network topologies. The results demonstrate
          the significant impact of the network topologies on the
          expansion plans and costs, as well as the computational
          benefits of the proposed approach.
	}%
	
	%\KEYWORDS{deterministic inventory theory; infinite linear programming duality; 
	%  existence of optimal policies; semi-Markov decision process; cyclic schedule}
	
	
	\maketitle
	%%%%%%%%%%%%%%%%%%%%%%%%%%%%%%%%%%%%%%%%%%%%%%%%%%%%%%%%%%%%%%%%%%%%%%
	
	\input{Introduction}

	\input{Formulation}
	\input{Relaxation}

	\input{SBD}
	\input{CG}

	\input{ModelDescription}

	\input{CaseStudy}

	\input{Conclusions}

	% References here (outcomment the appropriate case) 	
	% CASE 1: BiBTeX used to constantly update the references 
	%   (while the paper is being written).
	% \bibliographystyle{informs2014} % outcomment this and next line in 
	%	\bibliography{references.bib}

	%\bibliography{<your bib file(s)>} % if more than one, comma separated
	
	% CASE 2: BiBTeX used to generate mypaper.bbl (to be further fine tuned)
	%\input{mypaper.bbl} % outcomment this line in Case 2

	%\appendix
        % \input{Sections/Appendix}
	
\end{document}

%% file: Introduction.tex
%How do we derive a distribution for damage level. 

\section{Introduction}
%Within modern society, electric power plays a critically important role in ensuring the proper function of today's economic and social systems. %In spite of this role, extreme events, such as those posed by natural hazards like hurricanes, ice storms, and earthquakes \cite{wang_research_2016}, have the potential to severely disrupt the ability of power systems to deliver electricity.

The last decades have highlighted the vulnerability of the current
electric power system to weather-related extreme events. Between 2007
and 2016, outages caused by natural hazards, such as thunderstorms,
tornadoes, and hurricanes, amounted to 90 percent of major electric
disturbances, each affecting at least 50,000 customers (derived from
Form OE-417 of U.S. DOE). It is also estimated that 90 percent of all
outages occur along distribution systems \citep{Executive2013}.
Moreover, the number of weather-related outages is expected to rise as
climate change increases the frequency and intensity of extreme
weather events \citep{Executive2013}. Accordingly, it is critical to
understand how to harden and modernize distribution grids to
prepare for potential natural disasters.

Distributed Generation (DG) is one of the advanced technologies that
can be utilized to enhance grid resilience. DG refers to electric
power generation and storage performed by a collection of distributed
energy resources (DER). DG decentralizes the electric power
distribution by supplying power to the loads closer to where it is
located. The potential of DGs is realized via a system approach that
views DGs and associated loads as a microgrid
\citep{lasseter2002certs}. A microgrid is often defined as a
small-scale power system on medium- or low- voltage distribution
feeder that includes loads and DG units, together with an appropriate
management and control scheme supported by a communication
infrastructure \citep{resende2011service}. When faults occur in the
main grid, microgrids can be detached from the main grid and act in
island mode to serve critical loads by utilizing local DGs or work in
the grid-connected mode to provide ancillary services for the bulk
system restoration \citep{wang2016research}. Remotely controlled
switches (RCS), another advanced technology, can be used to increase
the grid flexibility by controlling the grid topology through a
communication network and facilitate microgrid formations in
emergencies. Other than the aforementioned operational enhancement
measures, a grid can also be hardened physically by installing
underground cables and/or upgrading the overhead lines with stronger
materials, which reduces the physical impact of catastrophic events
\citep{panteli2017}.

A critical issue in building resilient distribution grids is to
determine where to place such advanced devices (i.e., DGs, RCSs, and
underground cables) and which existing lines to harden.  It is also
important to understand the dependency between the distribution grid
and its associated communication network, which is critical to the
effective operation of a modernized grid during emergency situations
and is also vulnerable to extreme events
\citep{falahati2012reliability, gholami2016front,
  martins2017interdependence, li2017networked}.

To address this pivotal and pressing issue, this paper introduces the
Optimal Resilient Design Problem for Distribution and Communication Systems
(ORPDDC). The ORDPDC determines how to harden and
modernize an interdependent network to ensure its resilience against
extreme weather events. Like recent papers (e.g.,
\citet{yamangil2015resilient, Barnes2017Tools}, the ORDPDC 
takes into account a set of disaster scenarios, each defining a set of
power system components that are damaged during an extreme
event. These scenarios are generated from historical data or
probabilistic models of how power system components respond to
hazard-specific stress (e.g., wind speed and flood depth). The ORDPDC
considers the following upgrade options: a set of hardening
options on existing power lines and communication links and a set of
new components that can be added to the system---new lines, new
communication pathways, remotely controlled switches, and distributed
generation. The objective of the ORDPDC is to find the cheapest set of
upgrade options that can be placed on the grid in order to guarantee that a
minimal amount of critical and non-critical load be served in each
scenario. These guarantees are called the reslience criteria.

The ORDPDC is modeled with a two-stage stochastic mixed integer
program. The first stage decides an upgrade profile and the second
stage decides how to utilize the DGs, RCSs, and power
lines/communication links, whose availability is decided in the
first-stage, to restore critical loads up to resiliency criteria
(e.g., 98 \%) in each disaster scenario. For each scenario, the second
stage is viewed as a restoration model that identifies how to
reconfigure the grid. Within this second stage problem, the physics of
power flows is modeled with the steady-state, unbalanced three phase
AC power equations and constraints that ensure that the radial structure of
distribution grids is maintained. When the grid is reconfigured due to
some disturbances, each island or microgrid must be connected to at
least one control center that coordinates its DGs and loads and
operates its RCSs. This communication requirement is modeled with a
single-commodity flow model.

Several solution methods can be used to solve the ORDPDC,
taking advantage of its block diagonal structure. 
\citet{yamangil2015resilient} proposed a Scenario-Based
Decomposition (SBD) that restricts attention to a smaller set of
scenarios and adds new ones on an as needed basis (see Section
\ref{sec:sbd}). However, in the worst case, SBD must solve the
large-scale ORDPDC as a whole. Branch and Price (B\&P),
which combines column generation and branch-and-bound, is another
solution method for approaching large-scale mixed-integer programming
\citep{lubbecke2005selected}. Although widely successful on many
applications, it may suffer from degeneracy and long-tail effects as
problems become larger. To address these difficulties, several
stabilization techniques have been proposed and proven to be effective
in many applications (e.g., \citep{du1999stabilized,
  oukil2007stabilized, amor2009choice}). Nevertheless, the high degree
of degeneracy and the significant scale of the ORDPDC create
significant challenges for dual stabilization techniques.

To address these computational challenges, this paper proposes a B\&P
algorithm that systematically exploits the structure of the
ORDPDC. The algorithm starts with a compact reformulation that results
in strong lower bounds on the test cases and pricing subproblems that
are naturally solved in parallel. Moreover, the B\&P algorithm
tackles the degeneracy inherent in the ORDPDC through a variety of
acceleration schemes for the pricing subproblems: A pessimistic
reduced cost, an optimality cut, and a lexicographic objective. The
resulting B\&P algorithm produces significant computational
improvements compared to existing approaches.

The key contributions of this paper can be summarized as follows:

\begin{enumitemize}%[leftmargin=5mm]

\item The paper proposes the first planning model for resilient distribution
  networks that combines the use of advanced technologies (e.g., DGs,
  RCSs, and undergrounding) with traditional hardening options and
  captures the dependencies between the distribution grid and its
  associated communication system.

\item The paper proposes an exact B\&P algorithm for solving the ORDPDC
  problem, which systematically exploits the ORDPDC structure to
  obtain strong lower bounds and address its significant degeneracy
  issues.

\item The paper evaluates the impact of grid and communication system
  topologies on potential expansion plans. It also reports extensive
  computational results demonstrating the benefits of the proposed
  B\&P algorithm on the test cases.
\end{enumitemize}

The remainder of this paper is organized as follows. Section
\ref{sec:review} reviews related work on the ORDPDC. Section
\ref{sec:formulation} formalizes the ORDPDC and Section
\ref{sec:approxi} presents a tight linear approximation.  Section
\ref{sec:sbd} briefly reviews the SBD algorithm. Section \ref{sec:cg}
presents the new B\&P algorithm.  Section \ref{sec:data} describes the
test cases. Lastly, Section \ref{sec:casestudy} analyzes the behavior
of the model on the case studies and Section \ref{sec:performance}
reports on the computational performance of the proposed algorithm.
Section \ref{sec:conclusions} concludes the paper.

\section{Literature Review}
\label{sec:review}

There has been a considerable progress in advancing methods that
address weather-related issues at distribution level
\citep{wang2016research}. Many studies develop post-fault distribution
system restoration (DSR) models to bring power back as soon as
possible and restore critical loads after a severe outage. Recently,
DGs, RCSs, and redundant lines were utilized to leverage microgrids in
load restoration. Most of the studies assume the existence of those
devices beforehand \citep{chen2016resilient, ding2017new,
  gao_resilience-oriented_2016,yuan2017modified}. \citet{wang2016networked}
proposed a DSR model that utilizes the placement of dispatchable DGs.
The above-mentioned studies however propose post-contingency models. To
facilitate these novel restoration methods, the devices should be
placed in suitable places in advance. This paper focuses on the
optimal placement of those devices so that the grid survives potential
weather-related events.

Only a limited number of studies have discussed how to optimally add
resilience to distribution networks. Most relevant is the work by
\citet{Barnes2017Tools} and \citet{yamangil2015resilient} who propose
multi-scenario models for making a distribution grid resilient with
respect to a set of potential disaster scenarios. They propose
decomposition-based exact and heuristic solution approaches. However,
theses studies do not consider some of the upgrade options discussed
in this paper, and communication networks are not taken into
account. \citet{yuan2016robust} proposed a two-stage robust
optimization model by utilizing a bi-level network interdiction model
that identifies the critical components to upgrade for the resilience
against the $N-K$ contingency criterion. However, as pointed out in
\citet{Barnes2017Tools}, in practice, the computational complexity of
this approach grows quickly with the number of allowable faults. The
study also did not explicitly consider the dependency on the
communication network: A DG can supply power to the node it is placed
on and its children if they are not damaged by the attack.
\citet{carvalho2005decomposition} and \citet{xu2016placement} discuss
how to place RCSs in distribution systems, but only single fault
scenarios are assumed, which is not suitable for capturing
weather-related extreme events.

%Outside of the context of the resiliency, optimal placement of DGs in distirubtion systems has been actively studied in the last two decades [Review Prakash]. However, the studies assumed uncertainties arised from its generation output, load or electricity prices and the placement is not designed to be robust to the natural disaster-related disturbances [Wang 2014].

As the instrumentation of the grid increases, frameworks for modeling its
dependence on communication networks from a resilience viewpoint have
been studied \citep{martins2017interdependence,
  parhizi2015state}. \citet{resende2011service} proposed a
hierarchical control system, which assumes the existence of a controller
in each microgrid to allow for the coordination among distributed
generation units in the microgrid, while multiple microgrids are
organized by a central management controller. On the other hand,
distributed control systems are applied to microgrids where there are
many devices with their own controllers. Accordingly,
\citet{chen2016resilient} assumed that RCSs have local communication
capabilities to exchange information with neighboring switches over
short-range low-cost wireless networks and proposed a global
information discovery scheme to get the input parameters for a DSR
model. However, the assumption that RCSs are installed in all lines is
premature for current distribution systems.  \citet{wang2016networked}
proposed a two-layered communication framework where the lower-layer
cyber network supports microgrids where local control systems are
installed, while the upper-layer network is composed of multiple local
control systems that only communicate with their neighboring
counterparts. The study can be viewed as a hybrid of centralized and
decentralized framework: At a microgrid level, it is operated in a
centralized fashion, while the upper-level network is operated in a
decentralized manner. However, it did not consider fault scenarios in
communication networks. This paper only assumes the lower-layer cyber
network proposed in \citet{wang2016networked} by dynamically
allocating a local control system to each microgrid in islanding mode.
Moreover, this paper also considers potential faults in the
communication system.

To the best of our knowledge, this paper proposes, for the first time,
an exact optimization algorithm for expanding an integrated
distribution grid and communication network through the placement of
new DGs and RCSs and the hardening of existing lines in order to
ensure resilience against a collection of disaster scenarios.

%% file: Formulation.tex
\section{The ORDPDC}
\label{sec:formulation}

\begin{figure}[!t]
\begin{centering}
\includegraphics[width=0.50\paperwidth]{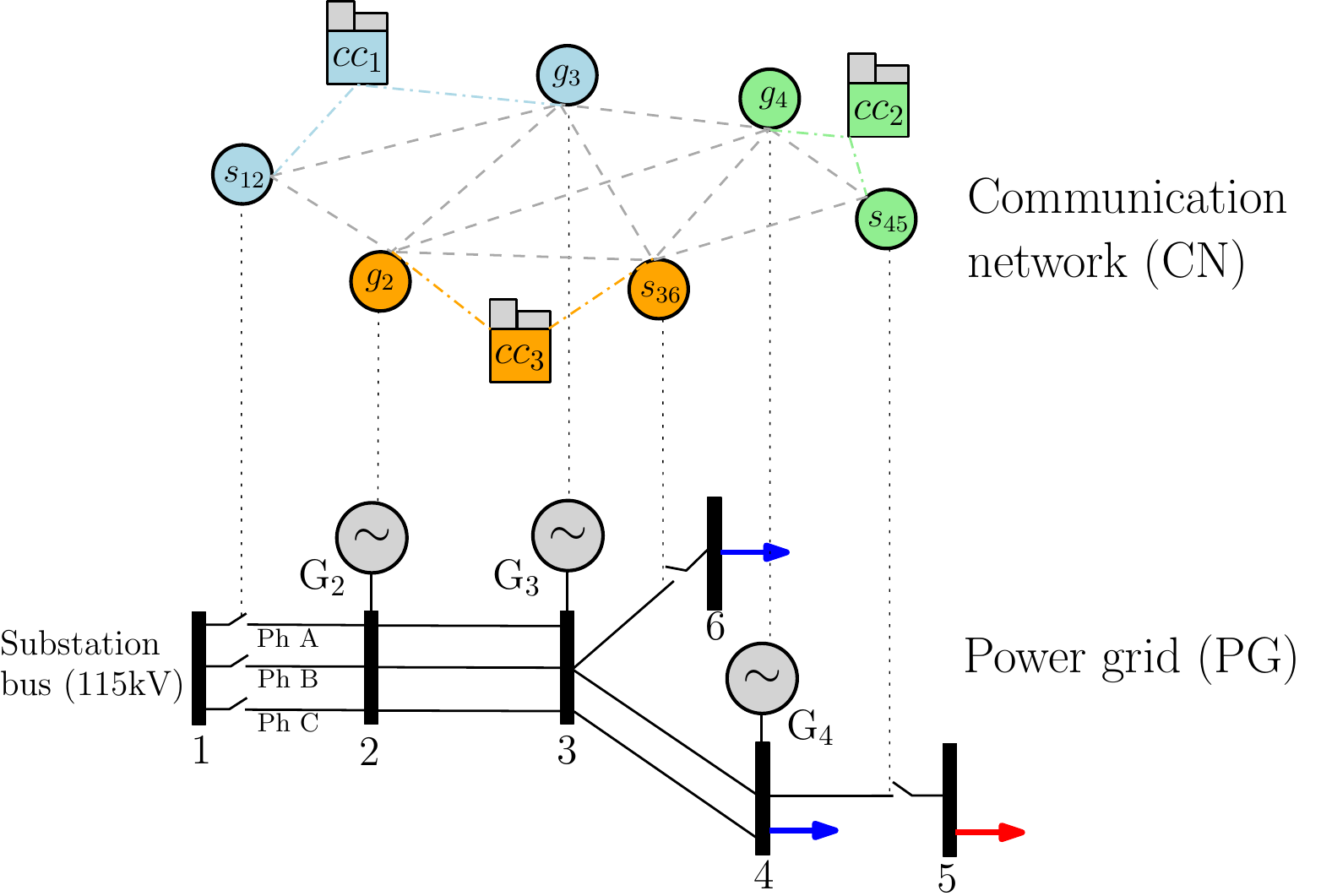}
\caption{The Cyber-Physical Network for Electricity Distribution. Solid lines represent power lines and dotted lines represent communication links.}
\label{fig:physical_cyber}
\end{centering}
\end{figure}

The ORDPDC considers an unbalanced three-phase distribution
grid coupled with a communication network, as illustrated in Figure
\ref{fig:physical_cyber}. In the figure, blue- and red-colored arrows
represent regular and critical loads. Nodes in the communication
networks may control a generator or a switch in the distribution
network, as indicated by dotted lines. The figure also
highlights how the line phases are interconnected at the buses and the
communication centers that will send instructions to generators and
switches remotely.
	
Let $G = (V,E)$ be an undirected graph that represents a distribution
grid and its available upgrade options: $V$ and $E$ denote the set of
buses and the set of distribution lines. The communication network,
along with its potential upgrade options, is represented by a
undirected graph $\Gc = (\Nc , \Ec )$, where $\Nc$ and $\Ec$ are the
set of communication nodes and a set of communication links. A
communication node is either a control point or an intermediate
point. Each control point is associated with some device in $G$ and
some nodes in $\Nc$ are designated as control centers.
	
The power grid $G$ depends on its communication network $\Gc$ in the
following way: A device in $G$ (e.g., a generator or a RCSs) is
operable only when its associated control point can receive a signal
from some control center in $\Gc$. This modeling enables islands to
form and to be operated independently only when at least one control
center can communicate to the island and, in particular, its
generator(s).

 %Define $\E_t^0 \subseteq E$ and $\E_t^n \subseteq E$ to be the set of lines in which a switch is installed and can be installed, respectively; define $\E_t := \E_t^0 \cup \E_t^n \subseteq E$. Likewise, define $\U^0$ and $\U^n$ to be the set of existing generators and generators can be installed, respectively; define $\U := \U^0 \cup \U^n$ and index it with $l \in \{1, \cdots, |U|\}$. Define $i: \U \rightarrow N$ to be a function that maps a generator to a bus where it is located. Let $\Nc_t \subseteq \Nc$ be a set of control points for switches and $\Nc_u \subseteq \Nc$ be that for generators; define $\Nc_c := \Nc_t \cup \Nc_u$ and $\tilde{i}: \E_t \cup \U \rightarrow \Nc_c$ to be the function that maps switch $e \in \E_t$ and generator $l \in \U$ to their control points in $\Gc$. 

Let $\G = (\N, \E)$ be the integrated system of $G$ and $\Gc$ with $\N
= N \cup \Nc$ and $\E = E \cup \Ec$. Let $\D$ be a set of damage
scenarios for $\G$ indexed with $\Sn := \{1, \cdots, |\D|\}.$ Each
scenario $s \in \Sn$ is a set of edges of $\E$ that are damaged under
$s$. The goal of the ORDPDC is to find an optimal upgrade
profile for the cyber-physical system $\G$ that is resilient with
respect to the damage scenarios in $\D$. The upgrade options include
a) the building of new edges in $\E$ (i.e., distribution lines or
communication links); b) the building of RCSs on some lines in $E$ to
provide operational flexibility; c) the hardening of existing edges in
$\E$ to lower the probability of damage, and d) the building of DGs at
some buses of the grid.
	
The ORDPDC is a two-stage mixed integer stochastic program. The
first-stage variables represent potential infrastructure enhancements
for the coupled network $\G$ and the second-stage variables capture
how upgrades serve the loads in each disaster scenario.
	
\subsection{Mathematical Formulation}

Table \ref{table:param} specifies the input data for the ORDPDC
problem, while Table \ref{table:var} describes the model
variables. The formulation assumes that all new lines come with
switches (i.e., $\E_x^n \subseteq \E_t^0$) which reflects current
industry practice.  Throughout this paper, an edge $e \in \E$ is
represented as an ordered pair $(e_h,e_t)$ for some $e_h, e_t \in \N$
and $\delta(e) = \{e_h, e_t\}$. The set of all edges incident to a
node $i \in \N$ is denoted by $\delta(i)$. The notation
$x_{\mathcal{A}}$ represents the projection of a vector $x$ to the
space of some index set $\mathcal A$, i.e., $(x_a)_{a \in
  \mathcal{A}}$: For instance, $x^s_{\E_x} = (x^s_e)_{e \in \E_x}.$
 
\begin{table}[!t]
\TABLE{The Parameters of the ORDPDC. \label{table:param}}
{\begin{tabular}{ll}
			\hline
			\up $G=(N, E)$ & an undirected extended distribution grid with available upgrade options \\%$N$ represents a set of buses and $E$ represents a set of lines.\\
			$\U := \U^0 \cup \U^n$ & a set of generators, indexed with $l$\\
			$ \qquad \U^0$  & a set of existing generators\\
			$ \qquad \U^n$   & a set of generators that can be installed\\
			$i(l) \in N$ & the bus in which the generator $l \in \U$ is located\\
			$\U_i \subseteq \U$ & the set of generators connected to bus $i \in N$\\
			$E_V \subseteq E$  & a set of transformers\\
			$\beta_e$ & maximum flow variation allowed between different phases on line $e \in E_V$\\
			$\C \subseteq 2^{|E|}$ & a collection of a set of edges which forms a cycle with a distinct node set \\
			$\Ph_e,\Ph_i,\Ph_l$ & a set of phases on line $e \in E$, bus $i \in N$, and generator $l \in \U$, respectively \\
			
			$T_{e}^k$ & a thermal limit on line $e \in E$ for phase $k \in \Ph_e$ \\
			$\underline V_{i}^k, \overline V_{i}^k$ & lower and upper bound on voltage magnitude at bus $i \in N$ on phase $k \in \Ph_i$ \\
			$Z_{e} = R_{e}+\mathbf{i} \ X_{e}$ & phase impedance matrix of line $e \in E$ \\
			$\mathcal{L} \subseteq N$ & a set of buses with critical loads \\
			$D_{i,p}^k + \mathbf{i}\ D_{i,q}^k$ & complex power demand at bus $i \in N$ on phase $k \in \mathcal P_i$ \\
			$\eta_{c}, \eta_{t}$ & resiliency criteria in percentage for critical and total loads respectively \\
			$\overline g^{k}_{l,p} + \mathbf{i} \ \overline g^{k}_{l,q}$ & complex power generation capacity of generator $l \in \U$ on phase $k \in \mathcal P_l$ \\
			
			$ \Gc=(\Nc, \Ec)$ & an extended associated communication network with potential upgrade options \\%$\Nc$ represents a set of communication nodes and $\Ec$ represents a set of communication links. \\
			$\Nc_c := \Nc_t \cup \Nc_u$ &\\
			$\qquad \Nc_t \subseteq \Nc$ & a set of control points for switches \\
			$\qquad\Nc_u \subseteq \Nc$ & a set of control points for generators\\
			$\tilde{i}(e) \in  \Nc_t, \tilde{i}(l) \in \Nc_u $ & the control point in $\Gc$ of a switch $e \in \E_t$ and a generator $l \in \U$, respectively\\
			$\tilde i_d \in \Nc$ & an artificial dummy node in $\Gc$\\

			$ \G=(\N, \E)$ & the integrated system of $G$ and $\Gc$\\%, where $\N = N \cup \Nc$, indexed with $i$, and $\E = E \cup \Ec$, indexed with $e$. \\
			
			$\E_x := \E_x^0 \cup \E_x^n$ & \\%a set of lines and links that are existing or can be installed.\\
			 $\qquad \E_x^0\subseteq \E$ & a set of existing lines and links\\
			 $\qquad \E_x^n\subseteq \E$ & a set of lines and links that can be installed \\
			$\E_t := \E^0_t \cup \E^{n}_t$ & \\%a set of lines in which a switch is installed, 
			$\qquad \E^0_t \subseteq E$ & a set of lines in which a switch is installed\\
			$\qquad\E^n_t\subseteq E$ & a set of lines in which a switch can be installed \\
			$\E_h \subseteq E$ & a set of lines or links that can be hardened \\
			$c^x_e$ & installation cost of $e \in \E_x^n$ \\
			$c^t_e$ & installation cost of switch on $e \in \E^n_t$ \\
			$c^h_e$ & line hardening cost of $e \in \E_h$ \\
			$c^u_l$ & installation cost of $l \in \U^n$ on the corresponding bus \\
%$e_h,\ e_t$ & the head and tail of $e \in \E$.\\
			\down $\D$ & a collection of sets of damaged lines for each scenario, indexed with $\Sn:=\{1, \cdots, |\D|\}$  \\%i.e., $\mathcal{D}_s$ corresponds to a subset of $\E$ consisting of damaged lines under scenario $s \in \Sn$ \\
					\hline
	\end{tabular}}{}
	\end{table}
	
	\begin{table}%[ht!]\fontsize{10}{10}\selectfont
		\TABLE{The Variables of the ORDPDC. \label{table:var}}
		{\begin{tabular}{ll}
			\hline
			\multicolumn{2}{l}{\up\down  \textbf{Binary variables}} \\
			 $x_e$ & 1 if $e \in \E_x^n$ is built\\
			$t_e$ & 1 if a switch is built on $e \in \E^n_t$ \\
			$h_e$ & 1 if $e \in \E_h$ is hardened \\
			$u_l$ & 1 if a generator $l \in \U^n$ is built. \\
			\multicolumn{2}{l}{\up For each disaster scenario $s \in \Sn$,} \\
			$z^s_e$ & 1 if $e \in \E$ is active during $s$ \\
			$x^s_e$ & 1 if $e \in \E_x$ exists during $s$ \\
			$t^s_e$ & 1 if a switch on $e$ is used or not  during $s$ \\
			$h^s_e$ & 1 if $e \in \E_h$ is hardened during $s$ \\
			$u^s_l$ & 1 if a generator $l \in \U^n$ is available during $s$\\
			$y^s_e$ & 1 if $i, \ j \in N$ can be disconnected, for $e = (i,j) \in C, \ C \in \mathcal{C}$, during $s$\\
			$b_e$ & 1 if the real power on line $e = (i,j) \in E$ flows from $j$ to $i$ during $s$\\
			$b'_e$ & 1 if the reactive power on line $e = (i,j) \in E$ flows from $j$ to $i$ during $s$\\
			\multicolumn{2}{l}{\up  \down \textbf{Continuous variables}} \\
			\multicolumn{2}{l}{For each disaster scenario $s \in \Sn$,} \\
			$d^{s,k}_{i} = d^{s,k}_{i,p} + \mathbf{i} \ d^{s,k}_{i,q}$ &  amount of power delivered at bus $i \in N$ on phase $k \in \Ph_i$ during $s$ \\
			$g^{s,k}_{l} = g^{s,k}_{l,p} + \mathbf{i} \ g^{s,k}_{l,q}$ &  amount of power generation of $l \in \mathcal{U}$ on phase $k \in \Ph_l$ during $s$ \\
			$s^{s,k}_{e,i} = p^{s,k}_{e,i} + \mathbf{i} \ q^{s,k}_{e,i}$ & power flow on $i$-end of line $ e \in E$, where $i \in \delta(e)$, on phase $k \in \Ph_e$ during $s$ \\
			$V^{s,k}_{i}$ & complex voltage at bus $i \in N$ on phase $k \in \Ph_i$ during $s$ \\
            $I^{s,k}_e$ & complex current on line $e \in E$ on phase $k \in \Ph_e$ during $s$ \\
			$v_{i}^{s,k}$ & squared voltage magnitude at bus $i \in N$ on phase $k \in \Ph_i$ during $s$ \\
			$f^{s}_{e}$ & the amount of artificial flow on $e \in \Ec$ during $s$\\
			\down$\gamma_{\tilde{i}}^s$ & indicator of connectivity of control point $\tilde i \in \Nc$ to some control center during $s$\\
			\hline
		\end{tabular}}{}
\end{table}

The presentation uses $w= ( x_{\E^n_x}, t_{\E^n_t}, h_{\E_h},
u_{\U^n})$ to denote upgrade profiles, $m$ the dimension of $w$, $c =
(c^x_{\E^n_x}, c^t_{\E^n_t}, c^h_{\E_h}, c^u_{\U^n}) \in \Real^m$ the
cost vector, and $w^s = ( x^s_{\E^n_x}, t^s_{\E^n_t}, h^s_{\E_h},
u^s_{\U^n})$ feasible upgrade profiles for each scenario $s \in
\Sn$. For each $s \in \Sn,$ $\Q(s)$ denotes the set of upgrade
profiles that enable the grid to maintain the predetermined load
satisfaction (resiliency) level $\eta_c, \eta_t$ (e.g., $\eta_c =
0.98$ and $\eta_t = 0.5$) under disaster scenario $s$.

With these notations, the ORDPDC is formulated as follows: 
% \todo[inline]{There is a subtlety here. The formulation $P$ assumes that $w \in \Q(s)$, given constraint 1b, but that is not always true in this case. Example, lets assume $w^1$ builds a line but $w^2$ does not. $w$ will be forced to build that line and $w$ could be infeasible for scenario 2 (maybe this introduces an unbreakable cycle in scenario 2).  I suggest two possible solutions 1. We explicitly state that we assume that $w \in \Q(s)$, given constraint 1b.  This assumption is ok here as we assume all new lines come with switches.  2. we add  the constraint $w \in \Q(s)$ to $P$ (and the later derivations of $P$)}
\begin{subequations}
\begin{alignat}{3}
(P) \quad &&\min \quad& c^T w \label{eq:first_obj}\\
	&&\mbox{s.t.}\quad & w \geq w^s, & \forall s \in \Sn, \label{eq:link_constr}\\
	&& 					&w^s \in \Q(s),  \ &\forall s \in \Sn, \label{eq:second_stage}\\
	&& 					&w \in \{0,1\}^{m}.\nonumber     
\end{alignat}
\label{eq:first}
\end{subequations}
\noindent
Problem ($P$) tries to find the optimal upgrade profile $w^* =
(x^*_{\E^n_x}, t^*_{\E^n_t}, h^*_{\E_h}, u^*_{\mathcal{U}^n})$ that
ensures resilient operations for each disaster scenario.  Equation
\eqref{eq:link_constr} ensures that an upgrade profile is feasible if
it dominates a feasible solution $w^s \in \Q(s)$ for each scenario
$s$, i.e., if the grid survives each of the extreme events in $\Sn$.
	
The set $\Q(s)$ is specified by resiliency constraints that are
expressed in terms of the AC power flow equations, load satisfaction
requirements, the communication network, and the grid topology:
\[
\Q(s) = \{ w^s \in \{0,1\}^m: (2),(3), \eqref{eq:generator/demand}, (5), \mbox{ and } \eqref{eq:topology}  \}
\]
where Constraints (2), (3), \eqref{eq:generator/demand}, (5), and
\eqref{eq:topology} are stated in detail in the following. The
variables in each $\Q(s)$ are indexed by $s$. For simplicity, this
section omits index $s$.
		
\subsubsection{Power Flow Constraints} 
\label{sec:formulation:powerflow}

\begin{figure}[!t]
\centering
\includegraphics[width = \textwidth] {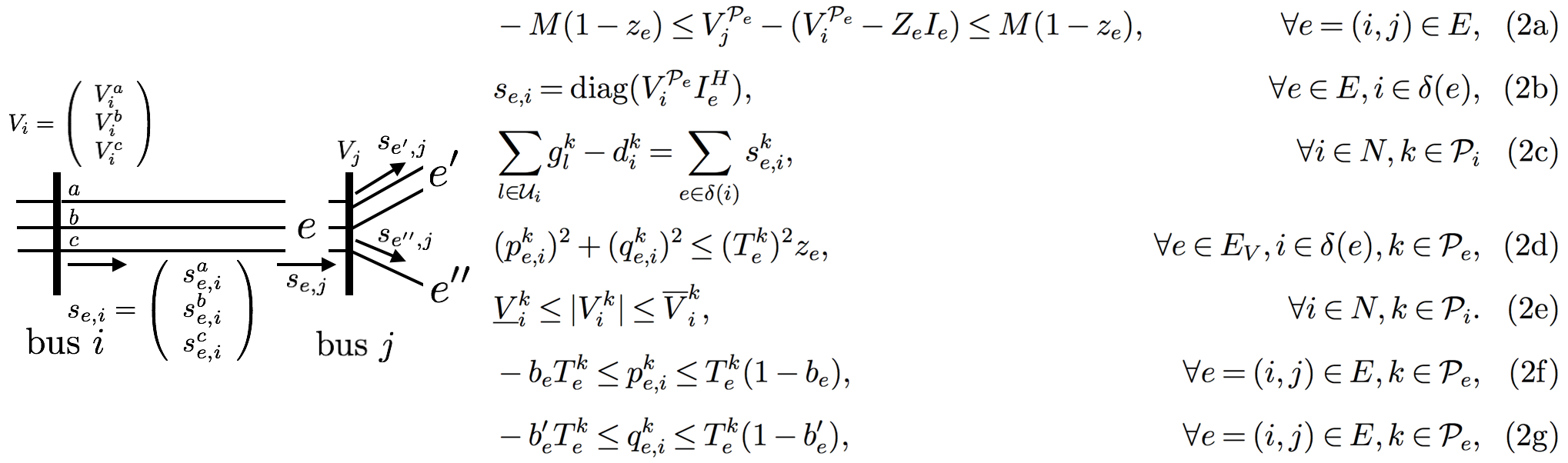}
\caption{Notations for the Power Flow Equations.}\label{fig:powerflow}
\end{figure} 

Figure \ref{fig:powerflow} specifies the power flow equations and
summarizes some of the notations. Let $\Ph = \{a, b, c\}$ denote the
three phases of the network. For each bus $i \in N$, define $V_i =
(V_i^k)_{k \in \Ph_i}$ and, for each line $e \in E$, define $I_e =
(I_e^k)_{k \in \Ph_e}$ and $s_{e,i} = (s_{e,i}^k)_{k \in \Ph_e}$. The
notations also use a superscript $\Ph' \subseteq \Ph$ to represent the
{\em projection} or the {\em extension} of a vector to the space of
$\Ph'$. For example, if $\Ph_i = \{a,b,c\}$ and $\Ph' = \{a,b\}$, then
$V_i^{\Ph'} = \left( V_i^a, V_i^b \right)^T$. If $\Ph_i = \{a,c\}$ and
$\Ph' = \{a,b, c\}$, then $V_i^{\Ph'} = \left( V_i^a, 0, V_i^c
\right)^T.$

% \begin{subequations}
% 	\begin{alignat}{2}
% 	%===============================================%
% 	%&\mathrm{\textbf{Power~flow~constraints}} \nonumber \\
% 	\hspace{2cm}&  -M(1-z_e) \le V^{\Ph_e}_j - ( V_i^{\Ph_e} - Z_{e} I_e) \le M(1-z_e), & \forall e = (i,j) \in E, \label{eq:Ohm}\\
% 	& s_{e, i} = \mbox{diag}(V_{i}^{\Ph_e} I_e^H),  & \forall e \in E, i \in \delta(e), \label{eq:Power}\\
% 	%			& -M (1-z_e) \le S_{e, j} - V_j^{\Ph_e} I_e^H  \le M (1-z_e),  & \forall e = (i,j) \in E, \label{eq:Power2}\\
% 	&\sum_{l \in \U_i} g^k_{l} - d^k_{i} = \sum_{e \in \delta(i)} s^k_{e, i},  & \forall i \in N, k \in \Ph_i \label{eq:power_balance}\\
% 	&  (p_{e,i}^{k})^2 +  (q_{e,i}^{k})^2 \le (T_e^k)^2 z_e,  &\   \forall e  \in E_V, i \in \delta(e), k \in \Ph_e, \label{eq:powerflow_thermal}\\
% 	%&  (p_{e,j}^{k})^2 +  (q_{e,j}^{k})^2 \le (T_e^k)^2 z_e,  \quad&\  \forall e = (i,j) \in E,  k \in \Ph_e, \label{eq:powerflow_thermal2}\\
% 	&\underline V_{i}^{k} \le |V_{i}^{k}| \le \overline V_{i}^{k}, &\  \forall i \in N, k \in \Ph_i. \label{eq:powerflow_voltage}\\	
% 	& -b_{e}T_{e}^{k}  \le p_{e,i}^{k} \le T_{e}^{k}( 1 - b_{e}), &  \forall e = (i,j) \in E,  k \in \Ph_e, \\
% 	& -b_{e}' T_{e}^{k}  \le q_{e,i}^{k} \le T_{e}^{k}( 1 - b'_{e}), &  \forall e = (i,j) \in E,  k \in \Ph_e, 
% 	\end{alignat}
% 	\label{eq:powerflow}%
% \end{subequations}

For each line $e = (i,j) \in E$, Ohm's law for 3-phase lines states
the relationship $V_j^{\Ph_e} = V_i^{\Ph_e} - Z_eI_e$ between $I_e$,
$V_i$, and $V_j$. For each line $e \in E$ and bus $i \in \delta(e)$,
the electric power flow equation $s_{e,i} =
\mbox{diag}(V_i^{\Ph_e}I_e^H)$ describes the relationship between
$s_{e,i}$, $V_i^{\Ph_e}$, and $I_e$, where superscript $H$ indicates
the conjugate transpose. In Figure \ref{fig:powerflow}, the big-M
method is used in Equations (2a) to apply Ohm's law only for available
lines; the big-M can be set as $\max_{j' \in \{i,j\}, k \in \Ph_e}
\overline{V}_{j'}^k - \min_{j' \in \{i,j\}, k \in
  \Ph_e}\underline{V}^k_{j'}$.  Equations (2c) is the balance equation
for power flow at each bus $i \in N$, i.e., the sum of incoming flows
equals the sum of the outgoing flows.

%Eq. (2a) states Ohm's law, and Eq. (2b) represents the electric power equation, where superscript $H$ means its conjugate transpose. 
%Define $S_{e,i}$ to be the complex power matrix of $i$-end of line $e  = (i,j) \in E$, where its diagonal entries correspond to the complex power flow $s_{e,i}^k$ at  $i$-end of the line for each phase $k \in \Ph_e$. The matrix is determined by the following electric power equation:
%\begin{equation}
%S_{e,i} = V_i^{\Ph_e}I_e^H,
%\label{eq:power_eq}
%\end{equation}
%where, $s_{e,i} = \mbox{diag}(V_i^{\Ph_e}I_e^H)$ and superscript $H$ means its Hermitian matrix. Let $p_{e,i} + \textbf{i} q_{e,i}$ be the retangular representation of $s_{e,i}$, where $p_{e,i} = (p_{e,i}^k)_{k \in \Ph_i}$ and $q_{e,i} = (q_{e,i}^k)_{k \in \Ph_i}$ denote the real and reactive power at $i$-end of line $e$. Moreover, we use the big-M method to force Eq. \eqref{eq:power_eq} to hold only for available lines:
%\begin{equation}
%-M (1-z_e) \le S_{e, i} - V_{i}^{\Ph_e} I_e^H  \le M (1-z_e), \quad \forall e \in E, i \in e,
%\end{equation}
%where $z_e$ represents the availablity of line $e \in E$.

%\begin{equation}
%\sum_{l \in \U_i} g^k_{l} - d^k_{i} = \sum_{e \in \delta(i)} s^k_{e, i}, \quad \forall  k \in \Ph_i, 
%\label{eq:power_balance}
%\end{equation}
%where $\delta(i)$ denotes the set of edges incident to bus $i \in N$.
%, $g_l^k = g_{l,p}^k + \textbf{i}  g_{l,q}^k$ denotes the complex power genration of generator $l$ on phase $k$, and $d_i^k = d_{i,p}^k + \textbf{i}  d_{i,q}^k$ represents the amount of complex power served at bus $i$ on phase $k$. 

Let $p_{e,i} + \textbf{i} q_{e,i}$ be the rectangular representation of
$s_{e,i}$, where $p_{e,i} = (p_{e,i}^k)_{k \in \Ph_i}$ and $q_{e,i} =
(q_{e,i}^k)_{k \in \Ph_i}$ denote the real and reactive power at the
$i$-end of line $e$. Equations (2d) and (2e) specify the thermal
limits on lines and the voltage bounds on buses.

%	\begingroup\makeatletter\def\f@size{10}\check@mathfonts

In some disaster scenarios when some of the lines are broken, power
flows of different phases on the same line can have opposite
directions, which is a very undesirable operationally. Equations (2f)
and (2g) prevent this behavior from happening.

The real and reactive power on different phase must stay within a
certain limit. Let $\widehat p_{e,i} = \sum_{\tilde{k} \in \Ph_e}
p_{e,i}^{\tilde{k}}$ and $\widehat q_{e,i} = \sum_{\tilde{k} \in
  \Ph_e} q_{e,i}^{\tilde{k}}$. Then, these limits are formulated as
follows:
\addtocounter{equation}{1}
\begin{subequations}
\begin{alignat}{2}
	%===============================================%
	%&\mathrm{\textbf{Power~flow~constraints}} \nonumber \\
	&\left(\underline{\beta}_{e}(1-b_{e}) + \overline{\beta}_{e} b_{e}\right)  \frac{\widehat p_{e,i}}{|\Ph_{e}|} \leq {p}_{e,i}^k \leq  \left(\underline{\beta}_{e}b_{e} + \overline{\beta}_{e}(1-b_{e})\right) \frac{\widehat p_{e,i}}{|\Ph_{e}|}, &\ \forall e  \in E_V, k \in \Ph_e, \label{eq:ph_bal_p}\\
	&\left(\underline{\beta}_{e} (1-b'_{e}) + \overline{\beta}_{e}b'_{e}\right)  \frac{\widehat q_{e,i}}{|\Ph_{e}|} \leq {q}_{e,i}^k \leq  \left(\underline{\beta}_{e}b'_{e} + \overline{\beta}_{e}(1-b'_{e})\right) \frac{\widehat q_{e,i}}{|\Ph_{e}|}, &\ \forall e  \in E_V, k \in \Ph_e, \label{eq:ph_bal_q}
	\end{alignat}
	\label{eq:ph_bal}%
\end{subequations}
where $\underline \beta_e = 1-\beta_e$ and $\overline \beta_e = 1+ \beta_2$.

\subsubsection{Generator/resiliency Constraints}  
\label{sec:formulation:gen}

Moreover, each generator $l \in \U$ has its own capacity and at least
some percentage of critical and total loads must be satisfied
as specified by the resiliency criteria $\eta_c$ and $\eta_t$.
\begin{subequations}
\label{eq:generator/demand}%
	\begin{alignat}{2}
	%&\mathrm{\textbf{Generator/resiliency constraints}} \nonumber \\
	& 0 \le g_{l,p}^k \leq \overline g^k_{l,p} u_l, \ g_{l,q}^k \leq \overline g^k_{l,q} u_l, & \forall l \in \mathcal{U}, k \in \Ph_l, \label{eq:capa_generator}\\
	&0 \le d_{i, p}^k \le D_{i,p}^{k}, \ 0 \le d_{i, q}^k \le D_{i,q}^{k}, &~\forall i \in N, k \in \Ph_i\label{eq:demand}\\
	& \sum_{i\in \mathcal{L}} d^{k}_{i,p} \geq \eta_{c} \sum_{i\in \mathcal{L}} D^k_{i,p},\quad \sum_{i\in \mathcal{L}} d^{k}_{i,q} \geq \eta_{c} \sum_{i\in \mathcal{L}} D^k_{i,q}, &\quad \forall k \in \Ph, \label{eq:critical_load_met}\\
	& \sum_{i\in N} d^{k}_{i,p} \geq \eta_{t} \sum_{i\in N} D^k_{i,p},\quad \sum_{i\in N} d^{k}_{i,q} \geq \eta_{t} \sum_{i\in N} D^k_{i,q}, & \forall k \in \Ph. \label{eq:ncritical_load_met}
	\end{alignat}
\end{subequations}
%\endgroup

Equation \eqref{eq:capa_generator} captures the power generation
capacity constraints. Equation \eqref{eq:demand} states that the
delivered power at each bus $i$ should not exceed the load. Equations
\eqref{eq:critical_load_met}-\eqref{eq:ncritical_load_met} enforce 
the resiliency constraints. 

\subsubsection{Communication Constraints} 
\label{sec:formulation:comm}

\begin{figure}[!t]
\begin{centering}
\includegraphics[width=\textwidth]{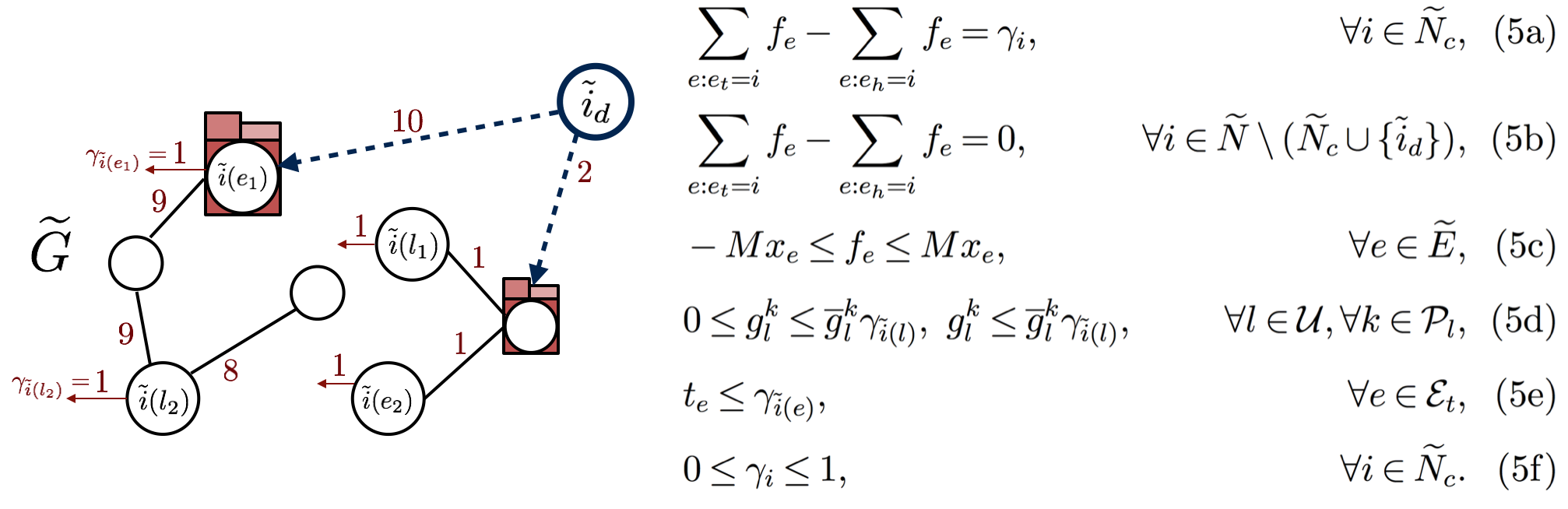}
\caption{The Single-Commodity Flow Model for $\Gc$ (red-colored squares denote control centers).}\label{fig:single_flow}
\end{centering}
\end{figure}

The operation of generators and RCSs depend on the communication
network: A generator $l \in \U$ and a RCS on line $e \in \E_t$ is
operable only if their associated control points $\tilde i(l) \in \Nc$
and $\tilde i (e) \in \Nc$ can receive a control signal from some
control centers through $\Gc$.  To capture the connectivity of a
vertex to some control centers, the formulation uses a
single-commodity flow model summarized in Equations (5) in Figure
\ref{fig:single_flow}. The formulation uses a dummy node $\tilde i_d$
to $\Nc$ and connect $\tilde i_d$ to all control centers with
additional links. The flow $f \in \mathbb{R}^{|\Ec|}$ originating from
the dummy node $\tilde i_d$ then is used to check the connectivity of
every node. By Equation (5c), the flow passes only through available
links during disaster $s$ (the big-M value is set to $|\Nc_c|$ in the
implementation). If a control point $i \in \Nc_c$ is connected with
some control center through some path, it can borrow a unit of flow
from $f$ to make $\gamma_{ i}$ 1, as specified in Equations (5a) and
(5b).  In other words, $\gamma_{ i}$ indicates whether control point
$i \in \Nc_c$ can receive a control signal.  If $\gamma_{ i}$ is 1,
the associated device in $G$ is operable by Equations (5d) and (5e).
 
%
%	\begin{subequations}
%	\begin{alignat}{2}
%	\hspace{6.4cm}&\sum_{e: e_t = i} f_{ e} - \sum_{e: e_h = i} f_{ e}= \gamma_{i}, & \ \forall  i \in \Nc_c, \label{eq:commu_flow_cp}\\
%	&\sum_{e: e_t = i} f_{ e} - \sum_{e: e_h = i} f_{ e}= 0, &\ \forall  i \in \Nc \setminus (\Nc_c \cup \{\tilde i_d\}),\label{eq:commu_flow}\\
%	&  -  M x_{ e} \le f_{ e} \le M x_{ e}, & \  \forall  e \in \Ec, \label{eq:artificialFlowCapa}\\
%	&0 \leq g^k_{l} \le \overline g^k_{l} \gamma_{\tilde{i}(l)},  \ g^k_{l} \le \overline g^k_{l} \gamma_{\tilde{i}(l)}, & \forall l \in \mathcal{U}, \forall k \in \Ph_l, \label{eq:commu_generator_0}\\
%	& t_e \le \gamma_{\tilde{i}(e)}, &\hspace{1cm} \forall e \in \E_t,\label{eq:commu_switch}\\
%	&0 \le \gamma_{i} \le 1, & \forall  i \in \Nc_c.
%	\end{alignat}
%	\label{eq:communication}%
%\end{subequations}

Some communication network may be affected by a failure in
distribution grid, e.g., when the grid fails to supply power to
communication centers. This kind of dependencies is not considered in
this paper but it can be easily captured if needed. Indeed, first
assign a small critical load to each communication center and add
constraints that restrict the auxiliary arcs between the dummy node
and each communication center to have positive flow only when the
associated communication center has a positive power supply. The
constraints can be expressed in terms of an extra binary variable for
each bus at which a communication center is located. The extra binary
variable determines if there is a positive power supply to the
communication center.

\subsubsection{Topological constraints.} 
\label{sec:formulation:topo}

The final set of constraints captures the topology restrictions in distribution systems:
\addtocounter{equation}{1}
\begin{subequations}
\begin{alignat}{2}
	& x_e \ge t_e, &\forall e \in \E,  \label{eq:topology_switch} \\
	& z_e = x_e - t_e, \hspace{4cm}& \forall e \in \E, \label{eq:line_availability}\\
	& x_e = h_e, &\forall e \in \D_s, \label{eq:topology_damage}\\
	& \sum_{e \in C} y_e  \leq |C|-1,&\forall C \in \mathcal{C},   \label{eq:topology_cycle1}\\
	& z_{\hat e} \le y_e,  & \forall \hat e \in E: \delta(\hat{e}) = \delta(e),  \ e \in C, \ C \in \mathcal{C}.  \label{eq:topology_cycle2}
\end{alignat}%
\label{eq:topology}%
\end{subequations}

Constraint \eqref{eq:topology_switch} restrict switches to be operable
only on existing lines. In Equation \eqref{eq:line_availability},
$z_e$ represents whether line $e \in \E$ is active under scenario
$s$. A line is active when it exists and its switch is off.  Equation
\eqref{eq:topology_damage} states that a damaged line during scenario
$s \in \Sn$ is inoperable unless it is hardened. Constraints
\eqref{eq:topology_cycle1} and \eqref{eq:topology_cycle2} ensures that
the distribution grid should operate in a radial manner. Accordingly,
Constraint \eqref{eq:topology_cycle1} eliminates the sub-tours within
$\C$. Since $G$ is usually sparse, the implementation enumerates all
the sub-tours $\C$ and variable $y_e$ indicates whether $i,j \in
\delta(e)$ are disconnected. If they are disconnected, then all the
lines between $i$ and $j$ are inactive by Constraint
\eqref{eq:topology_cycle2}.

Note also that, for existing lines not damaged under scenario $s$,
$x_e$ is fixed as one. For each line $e \in E \setminus \E_t$, $t_e$
is set to zero. Finally, for each line $e \in \E \setminus \E_h,$
$h_e$ is fixed as 0 and all the existing generators have $u_l = 1$.
This paper assumes perfect hardening, i.e., a hardened line survives
all disaster scenarios. This assumption can be naturally generalized
to imperfect hardening \citep{yamangil2015resilient}.

%% file: Relaxation.tex
\section{ Linearization of the ORDPDC}
\label{sec:approxi}

The formulation of the ORDPDC is nonlinear. This section discusses how to obtain an accurate linearization. 

\subsection{Linear Approximation of the AC Power Flow Equations for Radial Networks}
\label{sec:approxi:powerflow}

The main difficulty lies in linearizing constraints (2a--2b) for which
the formulation uses the tight linearization from
\citet{gan2014convex}.  The linearization is based on two assumptions:
(A1) line losses are small, i.e., $Z_{e}I_e I_e^H \approx 0$ for $e =
(i,j) \in E$ and (A2) voltages are nearly balanced, i.e., if $\Ph_i =
\{a, b, c\}$, then $V^a_i/V^b_i \approx V^b_i/V^c_i \approx
V^c_i/V^a_i \approx e^{i 2 \pi / 3}.$ Informally speaking, the
approximation generalizes the distflow equations to 3 phases, drops
the quadratic terms, and eliminates the current variables using the
balance assumption. The derivation assumes that all phases are
well-defined for simplicity. Moreover, if $A$ is an $n\times n$
matrix, then diag($A$) denotes the $n$-dimensional vector that
represents its diagonal entries. If $a$ is an $n$-dimensional vector,
then diag($a$) denotes the $n\times n$ matrix with $a$ in its diagonal
entries and zero for the off-diagonal entries. 

Let $s_i = \sum_{l \in \U_i} g_l - d_i$ denote the
power injection at bus $i$. By (A1), $s_{e,i} = s_{e,j}$ for
all $e \in (i,j) \in E$ and therefore, given $s_i$, $s_{e,i}$ ($i \in
\delta(e))$ is uniquely determined by Equation (2c). 

Now define $S_{e,i} := V_i I_e^H$, whose diagonal entries are $s_{e,i}$. 
Multipling both sides of $V_j = V_i - Z_e I_e$ with their conjugate transposes 
gives 
\begin{equation}
V_j V_j^H = V_iV_i^H - S_{e,i} Z_e^H - Z_e S_{e,i}^H + Z_e I_e I_e^H Z_e^H.
\end{equation}
By assumption (A1), this becomes 
\begin{equation}
V_j V_j^H = V_iV_i^H - S_{e,i} Z_e^H - Z_e S_{e,i}^H
\end{equation}
and, by restricting attention to diagonal elements only, 
\begin{equation}
v_j = v_i - \mbox{diag}(S_{e,i} Z_e^H - Z_e S_{e,i}^H).
\label{eq:Ohm2}
\end{equation}
where $(v^k_i)_{k \in \Ph_i} = \mbox{diag}(V_i V_i^H)$ represents the
squared voltage magnitude at bus $i \in N$.

By (A2), we have $S_{e,i} \approx \gamma^{\Ph_e} \mbox{diag}(s_{e,i}),$
where
$$\gamma = \left[\begin{array}{ccc} 1 & \alpha^2 & \alpha \\
\alpha & 1 & \alpha^2 \\
\alpha^2 & \alpha &1 \end{array} \right] \mbox{ and } \alpha = e^{-i 2 \pi /3}.$$ 

As a result, Equation \eqref{eq:Ohm2} can now be simplified as
follows: for each line $e = (i,j) \in E$ and $k \in \Ph_e$,
\begin{equation}
v^k_{i} =  v^k_{j} - \sum_{k' \in \Ph_e} 2\left[(\alpha^{n_k - n_{k'}}R_{e})^{kk'} p_{e,i}^{k'} + (\alpha^{n_k - n_{k'}}X_{e})^{kk'} q_{e,i}^{k'}\right],
	\label{eq:lindist}
\end{equation}
where $n_a = 2, n_b = 1, n_c = 0$, $R_e + \textbf{i} X_e = Z_e$, and superscript $kk'$ of a matrix denotes its $(k,k')$-entry.

In summary, Ohm's law and the power flow equation in Constraints (2a) and (2b) are approximated by Eq. \eqref{eq:lindist} for all $e = (i,j) \in E$ and $k \in \Ph_e$
and the big-$M$ is set to $\max_{j' =i,j} (\overline V_{j', k})^2 - \min_{j' =i,j} (\underline V_{j', k})^2$, along with Equation (2c). Accordingly, Constraint (2e) is replaced by the following constraint:
\begin{displaymath}
(\underline V_i^k)^2\le v_i^k \le (\overline V_i^k)^2, \quad \forall i \in N, k \in \Ph_i.
\end{displaymath}

\subsection{Linearization of \eqref{eq:ph_bal_p}-\eqref{eq:ph_bal_q}} 
\label{sec:appendix:ph_bal}

Constraints \eqref{eq:ph_bal_p} and \eqref{eq:ph_bal_q} contain
products of a binary variable and a bounded real variable. These
constraints are linearized without loss of accuracy using McCormick
inequalities \cite{mccormick1976computability}.

\subsection{Piecewise-Linear Inner Approximation of Thermal Limits}
\label{sec:approxi:thermal}

\begin{figure}[!t]
\centering
\includegraphics[width = 0.3 \textwidth]{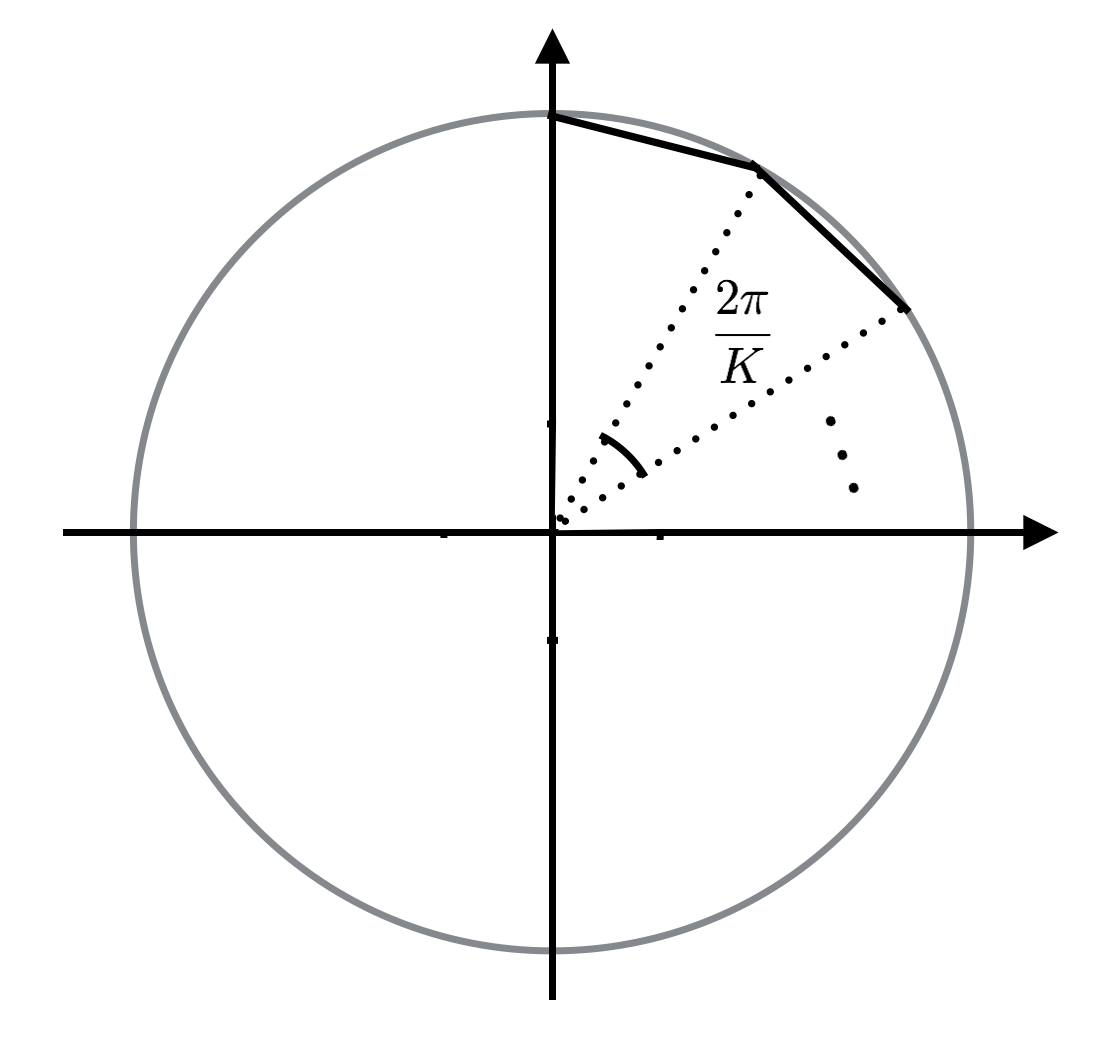}
\caption{The Piecewise-Linear Inner Approximation of a Circle.}
\label{fig:thermal}
\end{figure}

The quadratic thermal limit constraints (Constraint (2d)) can be approximated with $K$ linear inequalities as shown in Figure \ref{fig:thermal}. 
The resulting inequalities are as follows: for all $e \in E$, $i \in \delta(e)$, $k \in \Ph_e$:
\begin{subequations}	
\label{eq:linearlized_thermal}
	\begin{alignat}{2}
	&\left(\sin\left(\frac{2n\pi}{K}\right) - \sin\left(\frac{2(n-1)\pi}{K}\right)\right) p_{e,i}^k \nonumber\\
	&- \left(\cos\left(\frac{2n\pi}{K}\right) - \cos\left(\frac{2(n-1)\pi}{K}\right)\right) q_{e,i}^k \le  \sin\left(\frac{2\pi}{K}\right)T_{e,k}, &\  \forall n = 1, \cdots, K,\\
	&-M z^s_e \le p_{e,i}^{k} \le M z^s_e,\quad -M z^s_e \le q_{e,i}^{k} \le M z^s_e, & \quad \forall e \in E, k \in \Ph_e. \label{eq:lin_ther_i}
	%&-M z^s_e \le q_{e,j}^{k} \le M z^s_e,  -M z^s_e \le p_{e,j}^{k} \le M z^s_e,\quad \forall e \in E, k \in \Ph_e, \label{eq:lin_ther_j}
	\end{alignat}
\end{subequations}
%\endgroup
where the big-M is set to $\sum_{i \in N} D^p_{i,k}$. Our implementation
uses $K = 28$.

%% file: SBD.tex
\section{Scenario-Based Decomposition}
\label{sec:sbd}

In Section \ref{sec:performance}, the branch and price algorithm
presented in the next section is compared to the Scenario-Based
Decomposition (SBD) algorithm proposed by
\citet{nagarajan2016optimal}. SBD iteratively solves a master problem
$P(\Sn')$ which only includes the constraints of a subset of scenarios
$\Sn' \subseteq \Sn$. The algorithm terminates when the optimal
solution to $P(\Sn')$ is feasible (and hence optimal) for the
remaining scenarios $\Sn \setminus \Sn'$. Otherwise, at least one
scenario $s \in \Sn \setminus \Sn'$ is infeasible. Scenario $s$ is
added to $\Sn'$ and the process is repeated.

	%, which is described in Algorithm \ref{algo:SBD}. 
%	
%	\begin{algorithm}
%		\SetKwInOut{Input}{Input}
%		\Input{A set of disaster scenarios $\D$ indexed with $\Sn = \{1, \cdots, \CardiD\}$ and let $\Sn' = \{1\}$\;}
%		\Do{$s' \neq$ null} {$\sigma^* \leftarrow$ Solve $P(\Sn')$\;
%				Find $s' \in \Sn \setminus \Sn' $ in which $\sigma^*$ is infeasible\;
%				$\Sn' = \Sn' \cup \{s'\}$\;
%		}
%		\Return $\sigma^*$
%		\caption{Scenario Based Decomposition} \label{algo:SBD}
%	\end{algorithm}
%	The algorithm iteratively solves a master problem that consists of a subset of constraints for some scenarios $\Sn' \subseteq \Sn$, to be defined as $P(\Sn')$: 
%	\begin{alignat*}{3}
%	\left(P(\Sn')\right) \qquad && \min \quad& c^T w \\
%	&&\mbox{s.t.}\quad & w^s \le w, & \forall s \in \Sn',\\
%	&& &w^s \in \mathcal{Q}(s),  \ &\forall s \in \Sn',\\
%	&& &w \in \{0,1\}^{m}.
%	\end{alignat*}

%% file: CG.tex
\section{The Branch-and-Price Algorithm}
\label{sec:cg}

This paper proposes a branch-and-price (B\&P) algorithm for the
ORDPDC. The B\&P exploits the special structure of the ORDPDC 
in several ways. First, it uses a compact reformulation that yields a
better lower bound than the LP relaxation. The reformulation also
makes it possible to use column generation and solve independent
pricing problems associated with each scenario in parallel. Finally,
several additional techniques are used to accelarate the column
generation significantly. Section \ref{sec:reformulation} presents the
problem reformulation and Section \ref{sec:cg:basic} briefly reviews
the basic column generation of the B\&P algorithm. Section
\ref{sec:accelerating} introduces several acceleration schemes.  The
implementation of the B\&P algorithm is presented in Section
\ref{sec:implementation}.
	
\subsection{The Problem Reformulation }
\label{sec:reformulation}

Letting $\widetilde{\Q}(s)$ be the linearization of $\Q(s)$, the problem
($P$) is rewritten as 
\begin{subequations}
\begin{alignat}{3}
(P) \quad & \mbox{min} &  c^T w \nonumber\\
          & \mbox{ s.t. } & w - w^s \ge 0, &\ \forall s \in \Sn, \label{eq:link}\\
                && w^s \in \widetilde{\Q}(s),  &\ \forall s \in \Sn, \label{eq:sub}\\
         	&& w^s \in \{0,1\}^m,  &\ \forall s \in \Sn. \label{eq:binary}
\end{alignat}
\label{eq:CG:original}%
\end{subequations}

\noindent
Without the linking constraint \eqref{eq:link}, ($P$) can be
decomposed into $|\Sn|$ independent problems, each of which has a
feasible region defined by
 $$
\mathcal{P}^s =  \left\{w^s \in \mathbb{R}^m: \eqref{eq:sub} \mbox{ and } \eqref{eq:binary}\right\}, \ \forall s \in \Sn.
$$
Observe that $\mathcal{P}^s$ is bounded and let $\mathcal{J}^s =
\{\hat w^s_j \in \mathbb{R}^m: \hat w^s_j \mbox{ is a vertex of }
\mbox{conv} (\mathcal{P}^s) \}$ be the set of all vertices of
conv($\mathcal{P}^s$). Letting $\mathcal{J} = \cup_s \mathcal{J}_s$,
consider the following problem:

\begin{subequations}
\begin{alignat}{6}
(\widetilde P) \quad &\mbox{min} \quad & c^T w \nonumber\\
                     &\mbox{s.t.}   & w  & \quad - \sum_{j \in \mathcal{J}^s} \lambda_j^s \hat w_j^s \geq 0, \ & &&\forall s \in \Sn,\label{eq:link2}\\
				& & &\quad  \sum_{j \in \mathcal{J}^s} \lambda_j^s = 1, & && \forall s \in \Sn,\label{eq:convex}\\
			%	& & w \in & \{0,1\}\\
				& & &\quad w \in \{0,1\}^{m},\label{eq:binary2}\\
				& & & \quad \lambda^s_j \ge 0, & &&  \forall j \in \mathcal{J}^s,  s \in \Sn. 
\end{alignat}	
\label{eq:CG:reformulation}%
\end{subequations}
	
\begin{theorem}
($P$) and ($\Pa$) are equivalent.
\end{theorem}
\textbf{Proof.} Since ($P$) and ($\Pa$) have the same objective
function, it suffices to show that ($P$) has an optimal solution that
is feasible to ($\Pa$) and vice versa. Let $(\bar{w},
\{\bar{w}^{s}\}_{s \in \Sn})$ be the optimal solution of ($P$). By the
Farkas-Minkowski-Weyl theorem \citep{schrijver1998theory}, $\bar w^s$
can be expressed as a convex combination of some extreme points in
$\mathcal J^s$, for each $s \in \Sn$. Hence, we can construct a
feasible solution of ($\Pa$) from $(\bar{w}, \{\bar{w}^{s}\}_{s \in
  \Sn})$.
	
Consider now an optimal solution of ($\Pa$), $(\bar{w}^\prime,
\{\bar{\lambda}^{s \prime}_j\}_{j \in \mathcal{J}^s} \mbox{ for } s
\in \Sn)$. By \eqref{eq:link2}, if $\bar{\lambda}^{s \prime}_j > 0$
for $j \in \mathcal{J}^s$, $\hat{w}^s_j$ is dominated by
$\bar{w}^\prime.$ Therefore, it is possible to construct another
optimal solution to ($\Pa$) by choosing a single $j^*$ for which
$\bar{\lambda}^{s \prime}_{j^*} > 0$ for each $s \in \Sn$, setting
$\bar{\lambda}^{s \prime}_{j^*}$ to one and setting the other
$\bar{\lambda}^{s \prime}_{j}$ to zero. By definition of
$\mathcal{J}_s$, the constructed optimal solution is feasible for
($P$). $\square$
	
\noindent
This paper uses a branch-and-price algorithm to solve ($\Pa$). Let
$LP_{\Pa}$ denote the LP relaxation of ($\Pa$). Since the feasible
region of ($\Pa$) is the intersection of the convex hulls of each
subproblem, $LP_{\Pa}$ yields a stronger lower bound than the LP
relaxation of ($P$).
	
\subsection{The Basic Branch and Price }
\label{sec:cg:basic}

The B\&P algorithm uses a restricted master problem ($M$) with a
subset of columns of ($\Pa$) and $|\Sn|$ independent subproblems
($P_s$) for $s \in \Sn$, instead of handling $LP_{\Pa}$ globally. The
column generation starts with an initial basis that consists of the
first-stage variables $w$, a column associated with a feasible
solution for each subproblem, and some slack variables. Let
$\widetilde{\mathcal{J}}^s$ be the corresponding subset of
$\mathcal{J}^s$. The restricted mater problem ($M$) is
as follows:
\begin{subequations}
\begin{alignat}{5}
(M) \quad &\mbox{min} \ & c^T w \nonumber\\
	&\mbox{s.t.} \  & w - & \sum_{j \in \widetilde{\mathcal{J}}^s} \lambda_j^s \hat w_j^s \geq 0, & \ \forall s \in \Sn,\label{eq:m_link}\\
			& & &\sum_{j \in \widetilde{\mathcal{J}}^s} \lambda_j^s = 1, &  \forall s \in \Sn,\label{eq:m_convex}\\
			%	& & w \in & \{0,1\}\\
			& && w \ge 0,\ \\
			& && \lambda^s_j \ge 0, \quad  &  \forall j \in \widetilde{\mathcal{J}}^s, \ \forall s \in \mathcal{S}.
\end{alignat} 
\label{eq:master}
\end{subequations}
\noindent
and the pricing problem for scenario $s$ is specified as follows:
\begin{equation*}
\begin{array}{lll}
(P_s) \quad & \  \min  & \ - \bar{\sigma}^s \ + {\bar{y}^{sT}} \ w^s \\
& \mbox{s.t.} & \ w^s \in \widetilde{\Q}(s), \\
&             & \ w^s \in \{0,1\}^m,
\end{array}
\end{equation*}

\noindent
where, for a scenario $s$, $\bar{y}^{s}$ is the dual solution for
constraints \eqref{eq:m_link} and $\bar{\sigma}^s$ is the dual
solution of the convexity constraint (\ref{eq:m_convex}).

\subsection{Acceleration Schemes } 
\label{sec:accelerating}

The performance of column generation deteriorates when the master
problem exhibits degeneracy, leading to multiple dual solutions which
may significantly influence the quality of columns generated by the
pricing problem. The master problem $(M)$ suffers from degeneracy,
especially early in the column-generation process. Initially, ($M$)
has $(m+1) |\Sn|$ constraints, $m$ columns corresponding to the
first-stage variables $w$, and $|\Sn|$ columns for the second-stage
variables $\{\lambda^s\}_{s \in \Sn}$. Therefore, in early iterations,
linear solvers have a natural tendency to select $m(|\Sn|-1)$ columns
from the slack variables in Constraints \eqref{eq:m_link}. For example, assume that
the slack variable is in basis for the constraint involving a non-basic first-stage variable $w_k$ and
a scenario $s$ in Constraints \eqref{eq:m_link}. By complementary
slackness, this implies that the dual variable is zero. Consider a
vertex $\hat w^s$ whose $k$-th entry is non-zero. The value $\bar
y^s_k w^s_k$ is zero in the pricing problem. However, for this vertex
to enter the basis, it must incur the cost $c_k$ of $w_k$, which is
ignored in the pricing subproblem.  As a result, subproblem ($P_s$)
prices many columns too optimistically and generates columns that do
not improve the current objective value, resulting in a large number
of iterations.
	  
\subsubsection{Pessimistic Reduced Cost }
\label{sec:accelerating:revised_rc}

In order to overcome the poor pricing of columns, this section first
proposes a pessimistic pricing scheme that selects more meaningful
columns in early iterations. Consider a solution $w^s$ to the pricing
problem. If $w^s_k = 1$ but the first-stage variable $w_k$ is not in
basis, then by the relevant constraint from $\eqref{eq:m_link}$, the variable
$\lambda_j^s$ corresponding to $w^s$ can only enter in the basis at 1
if $w_k$ is also in the basis at 1. As a result, the pessimistic
pricing scheme adds the reduced cost 
$
c_k - \sum_{s \in \Sn} \bar{y}^s_k 
$
to the pricing objective, which becomes 
\[
- \bar{\sigma}^s + (\bar{y}^s)^T w^s + \sum_{k \in \eta} (c_k - \sum_{s \in \Sn}\bar{y}^s_k) w^s_k \label{eq:revised_rc}
\]
where $\eta$ is the set of non-basic first-stage variables, i.e., $
\eta = \{ k \ | \ w_k \mbox{ is non-basic} \}. $ Note that column
generation with this pessimistic pricing subproblem is not guaranteed
to converge to the optimal linear relaxation. Hence, the
implementation switches to the standard pricing problem in later
iterations.

\subsubsection{Optimality Cut} 
\label{sec:accelerating:opt_cut}

A solution to the master problem ($M$) where the first-stage variables
take integer value gives an upper bound to the optimal solution. The
B\&P algorithm periodically solves the integer version of ($M$) to
obtain its objective value $\bar v(M)$. The constraint 
\[
c^T w^s \le \bar v(M)
\]
can then be added to the pricing subproblem for scenario $s$ since any
solution violating this constraint is necessarily suboptimal. As shown
later on, this optimal cut is critical to link the two phases of the
column generation, preventing many potential columns to be generated
in the second phase.

\subsubsection{A Lexicographic Objective for Pricing Subproblems} 
\label{sec:accelerating:multi_obj}

In general, sparse columns are more likely to enter the basis in the
master problem ($M$). As a result, the B\&P algorithm uses a
lexicographic objective in the pricing subproblem. First, it minimizes
the (pessimistic or standard) reduced cost. Then it maximizes sparsity
by minimizing $1^T w^s$ subject to the constraint that the reduced
cost must be equal to the optimal objective value of the first stage. 

%The over-optimistic pricing of a column is due to the objective
%function of many zero's for ($P_s$). For , however, a solution
%with least positive entries has a more possibility to enter the basis,
%i.e., improve the objective value. In order to get a sparse column, we
%add a variable $\varsigma^s$, which represents the sum of all the
%binary variables, to the subproblem $(P_s)$ with small enough
%objective coefficient $c_\varsigma$, e.g. 1e-2. We then subtract
%$c_\varsigma \bar{\varsigma}^s$ from the objective value to check if
%the subproblem has a negative reduced cost. We exclude the variable
%when a subproblem has a trivial objective function with all zero's.
%	
%After switching to the original reduced cost, if it takes a longer
%time, e.g. $> 20$ seconds, in solving ($P_s$) with the penalty term,
%we exclude the term and solve the problem in two steps. If a
%subproblem ($P_s$) has a negative reduced cost, we resolve ($P_s$)
%after including an additional constraint that the pricing value has to
%be less than or equal to 0 and replacing the objective function as
%minimizing $\varsigma^s$.

\subsection{The Final Branch and Price Implementation} 
\label{sec:implementation}

\subsubsection{Column Generation} 
\label{sec:implementation:cg}

The column generation starts with an initial basis built from the
optimal solutions of each subproblems under the objective function of
$c^T w^s$. It then proceeds with two phases of column generation,
first using the pessimistic reduced cost and then switching to the
standard one. 

The second phase terminates when the optimality gap becomes lower than
the predetermined tolerance, e.g., $0.1\%$. The lower bound is based
on Lagrangian relaxation. Given a pair $\bar{w}$ and $(\bar{y},
\bar{\sigma})$ of optimal primal and dual solutions for ($M$), the
Lagrangian relaxation is given by
\[
L(\bar{w},\bar{y},\bar{\sigma}) = c^T \bar{w} + \sum_{s \in \mathcal{S}} {\cal O}_s(\bar{y},\bar{\sigma})
\]
where ${\cal O}_s(\bar{y},\bar{\sigma})$ is the optimal solution of
the pricing problem for scenario $s$ under dual variables
$(\bar{y},\bar{\sigma})$. The first phase uses the same technique for
termination, although the resulting formula is no longer guaranteed to
be a lower bound. Once the gap between the upper bound and the
``approximate'' lower bound is smaller than the tolerance, the column
generation process moves to the second phase.

The column generation also avoids generating dominated columns. Assume
that $[w^s_1 = 1, w^s_2 = 1]$ is a feasible solution of $(P_s)$ and
the corresponding column has been added to the master problem ($M$).
Then, there is no need to consider a solution $[w^s_1 = 1, w^s_2 = 1,
w^s_3 =1]$. The column generation adds the constraint of $w^s_1 +
w^s_2 \le 1$ to $(P_s)$ when such a dominated solution is produced and
does not include it in the master problem.
	
%	\begin{align}
%	\label{eq:Lagrangian}
%	L(\bar{y}) %&= c^T \bar{w} + \sum_{s \in \mathcal{S}}  \mbox{min}_{j \in J^s} \left( - \bar{\sigma}^s - \bar{y}^{sT} \hat{w}_j^s \right)\\
%	& = c^T \bar{w} + \sum_{s \in \mathcal{S}} \mbox{ optimal value of $(P_s)$ under $(\bar{y}^T, \bar{\sigma})$}.
%	\end{align}
	
%	We use $L(\bar{y})$ as a lower bound on the optimal value of $LP_{\Pa}$. Note that it only involves a trivial calculation since it is the sum of the optimal objective value of the current restricted master problem ($M$) and that of subproblems ($P_s$)'s. 
%	
%	In addition, when the revised reduced cost (\ref{eq:revised_rc}) is applied, due to the optimality of $\bar{y}$ for the dual of the current ($M$), the additional term $\sum_{k \in \eta} \left(c_k - \bar{y}^T A_{\cdot k}\right) \hat{w}^s_{jk}$ must be nonnegative. Therefore, the sum of the optimal objective value of ($M$) and that of ($P_s$)'s, denoted by $L'(\bar{y})$, can also be used as a relaxed lower bound for the first phase of the algorithm at which \eqref{eq:revised_rc} is applied. Hence, we switch to the actual reduced costs when the gap between the current objective value and $L'(\bar{y})$ is less than 0.1\%.
%	
\subsubsection{The Branch and Bound} 
\label{sec:implementation:bnb}

After convergence of the column generation to $LP_{\Pa}$, the branch
and bound algorithm solves the restricted master problem ($M$) with the
integral condition $w \in \{0,1\}^m$ to obtain a strong primal bound.
In general, this incumbent solution is of very high quality and the
average optimality gap is 0.19\%. Therefore, the branch and price
algorithm uses a depth-first branch and bound. Moreover, at each branching node,
it selects the variable that minimizes the optimality gap.

%% file: ModelDescription.tex
\section{Description of the Data Sets} 
\label{sec:data}

This section describes the distribution test systems. The data set is
available from \url{https://github.com/lanl-ansi/micot/} in the
\url{application_data/lpnorm} directory.  Details of the data format
are available from
\url{https://github.com/lanl-ansi/micot/wiki/Resilient-Design-Executable}.

\begin{figure}[!t]
\centering
\subfloat[Urban]{\includegraphics[width=0.4\textwidth]{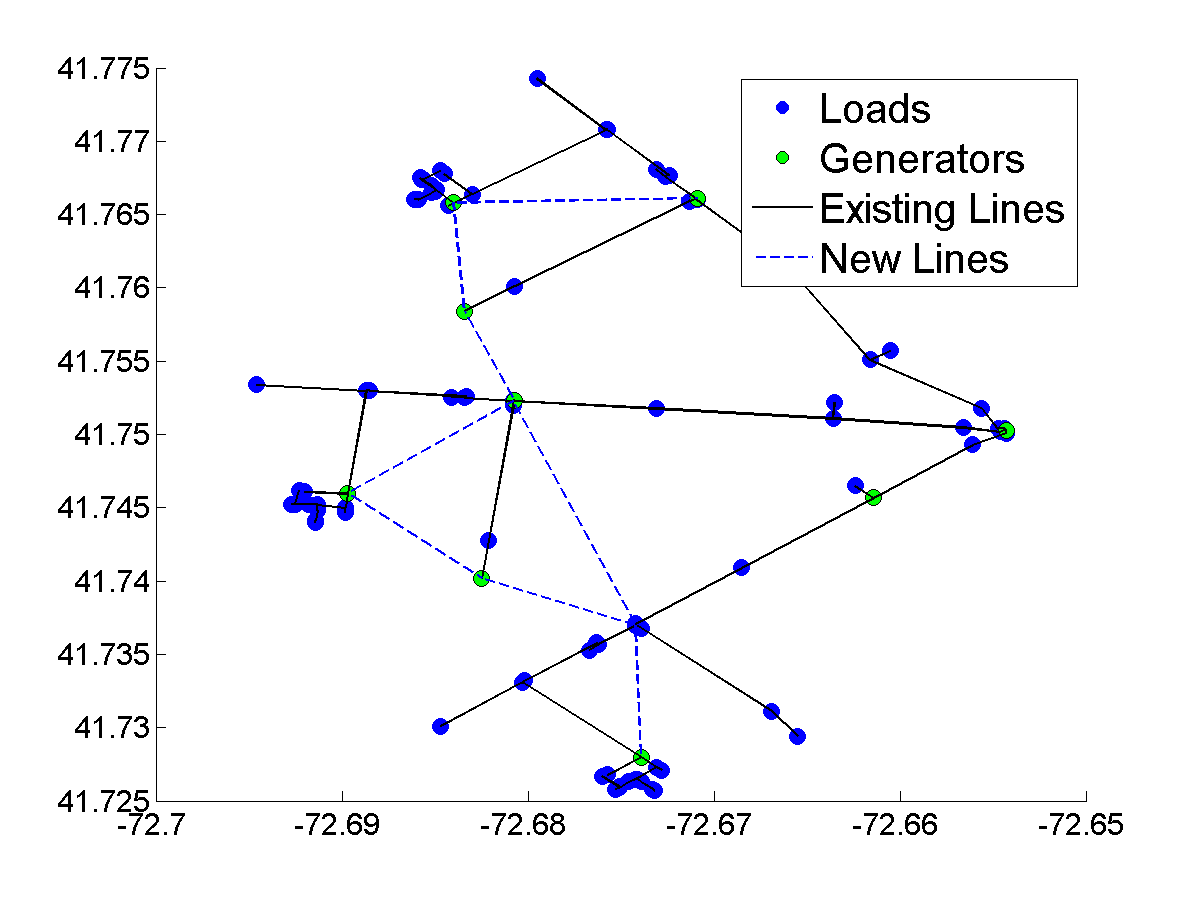}}\label{fig:1:a}
\subfloat[Rural]{\includegraphics[width=0.4\textwidth]{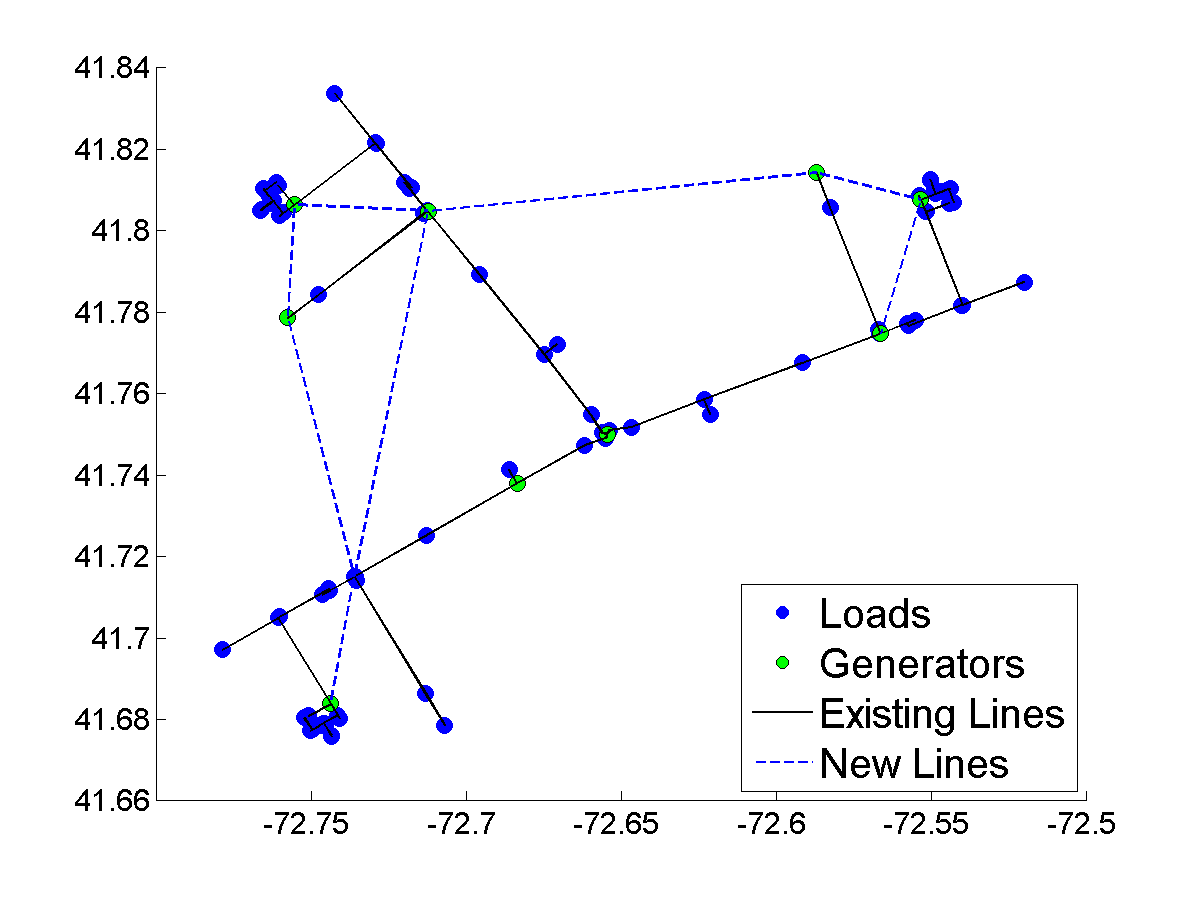}}\label{fig:1:b}
\caption{The urban and rural distribution systems which contain three copies of the IEEE 34 system to mimic situations where there are three normally independent distribution circuits that support each other during extreme events. These test cases include 109 nodes, 118 generators, 204 loads, and 148 edges.}
\label{fig:ieee34}
\end{figure}

The first two sets, the {\em Rural} and {\em Urban} systems, is from
\citet{yamangil2015resilient}. They are based on the IEEE 34 bus
system \citep{Kersting1991} (see Figure \ref{fig:ieee34}) and
replicate the 34-bus distribution feeder three times.  All three
feeders are connected to a single transmission bus and candidate new
lines were added to the network to allow back-feeds.  In the rural
model, the distribution feeder was geolocated to model feeders with
long distances between nodes. Similarly, the urban network was
geolocated to model compact feeders typical of urban
environments. Geolocation of these networks has the net effect of
adjusting the lengths of the power lines and their associated
impedance values.  Spreading the network out also increases the
hardening and new line costs. As a result, the rural system is
expected to favor solutions with distributed generation and the urban
system solutions with new lines and switches (in addition to hardening
lines). The fixed cost of installing a new distributed generator is
set at \$500k. The cost of a distributed generator is set at \$1,500k
per MW based on the 2025 projections from
\citet{U.S.EnergyInformationAdministration2014}. The cost of
installing new switches for 3-phase lines is set between 10k and 50k
\citep{Ba}. The cost of new underground 3-phase lines is set at about
\$500k per mile and the cost of new underground single phase lines is
set at about \$100k per mile. The hardening cost was set at roughly
\$50k and \$10k per mile for multi-phase and single-phase lines
\citep{Of2005}. The third network, \textsc{network123}, is based on
the 123-node network of \citet{Kersting1991}. This network was
unaltered except for adding new line candidates and labeling large
loads as critical.

The communication network $\Gc$ is built to conform to $G$. Let $G' =
(N',E')$ be the duplicate of $G$. For each generator $l \in \U$, its
duplicate $i(l)$ represents its control point. Consider $\E'_t
\subseteq E'$, the duplicate of $\E_t$. To represent the control point
for a switch, $e \in \E'_t$ is divided in the middle and a new vertex
$v_e$ is added to represent the control point for the switch.  In
other words, the edge $e = (e_h, e_t) \in \E'_t$ is replaced by a new
vertex $ v_e$ and two new edges $ e_1 = (e_h, v_e), e_2 = ( v_e,
e_t)$. The test cases assume that the damage, installation, and
hardening of a line in $G$ are also incurred for the corresponding
line in $\Gc$. These assumptions can be easily generalized without changing
the nature of the model. 

The experimental evaluation considers 100 scenarios per damage
intensity for all three networks and the damage intensities are taken
in the set $\{1\%,$ $2\%,$ $3\%,$ $4\%,$ $5\%,$ $10\%,$ $15\%,$
$20\%,$ $25\%,$ $30\%,$ $35\%,$ $40\%,$ $45\%,$ $50\%,$ $55\%,$
$60\%,$ $65\%,$ $70\%,$ $75\%,$ $80\%,$ $85\%,$ $90\%,$ $95\%,$
$100\%\}$. The scenario generation procedure is based on damage caused
by ice storms. The intensity tends to be homogeneous on the scale of
distribution systems \citep{Sa2002}. Ice storm intensity is modeled as
a per-mile damage probability, i.e. the probability at least one pole
fails in a one mile segment of power line. Each line is segmented into
1-mile segments and a scenario is generated by randomly failing each
segment with the specified probability. This probability is normalized
for any line segment shorter than 1 mile. A line is ``damaged'' if any
segment fails.

%% file: CaseStudy.tex
\section{ Case Study} 
\label{sec:casestudy}

This section analyzes the behavior of the optimization model on a
variety of test cases. In particular, it studies how the topology of
the distribution grid and the dispersion level of its communication
network affect the optimal design. For each network described in
Section \ref{sec:data}, this section analyzes the optimal design under
different settings of damage probability, the resiliency level, and
the number of communication centers. The default value of $\eta_c$ and
$\eta_t$ are 98\% and 50\% respectively, the default number of
communication centers is 4, and the phase variation parameter $\beta$
is set to 15\% for $E_V$ and $\infty$ otherwise. Unless specified otherwise, the comparisons are based
on these default values.

\subsection{Impact of grid topology} \label{sec:impact_topology:grid}

\begin{figure}[t!]
%	\begin{figure}[b!]
	\centering
	\subfloat[Rural network] {\includegraphics[width=0.46\textwidth]{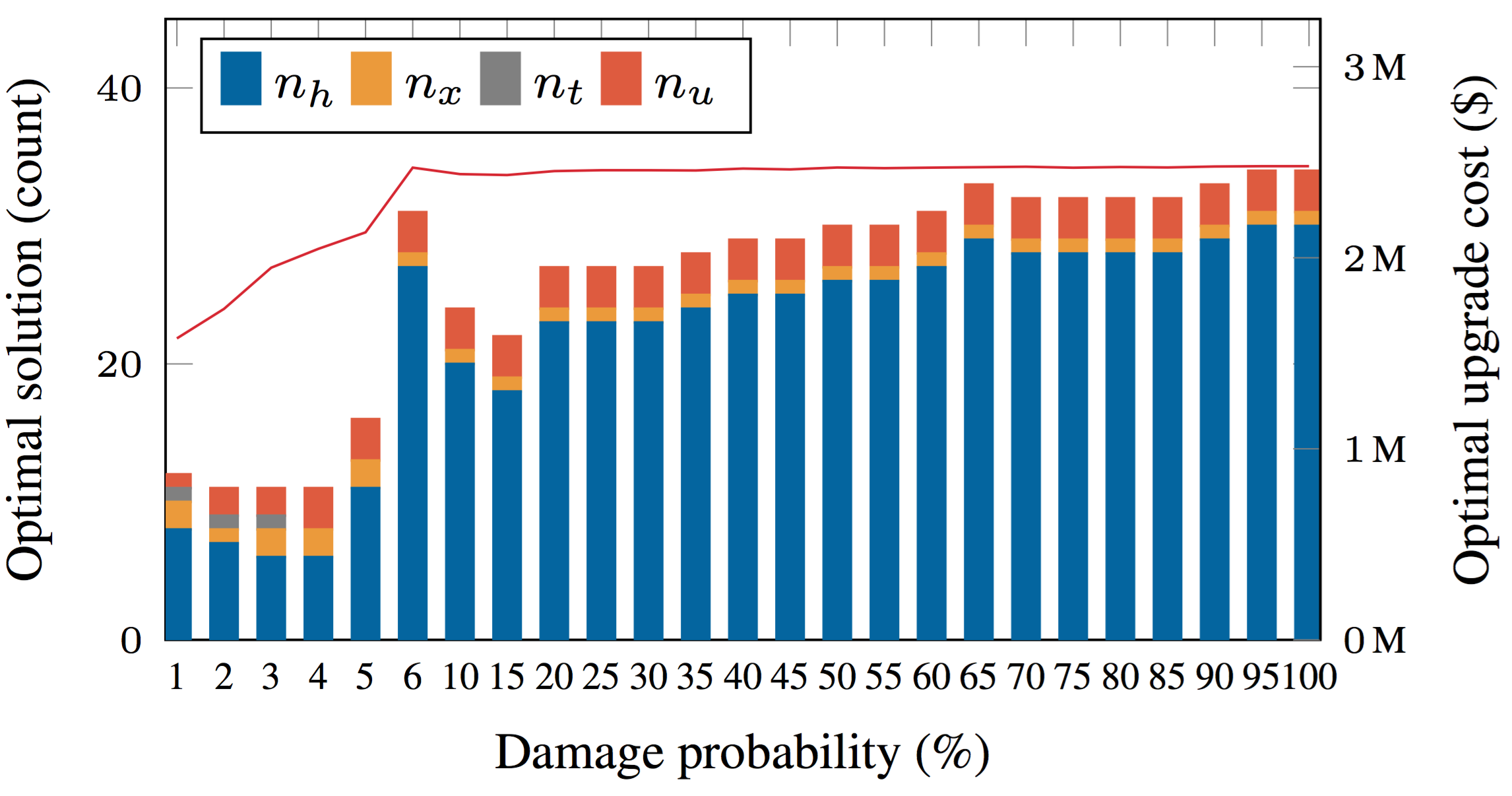}}\label{fig:opt_sol:rural}
	\subfloat[Urban network]{\includegraphics[width=0.46\textwidth]{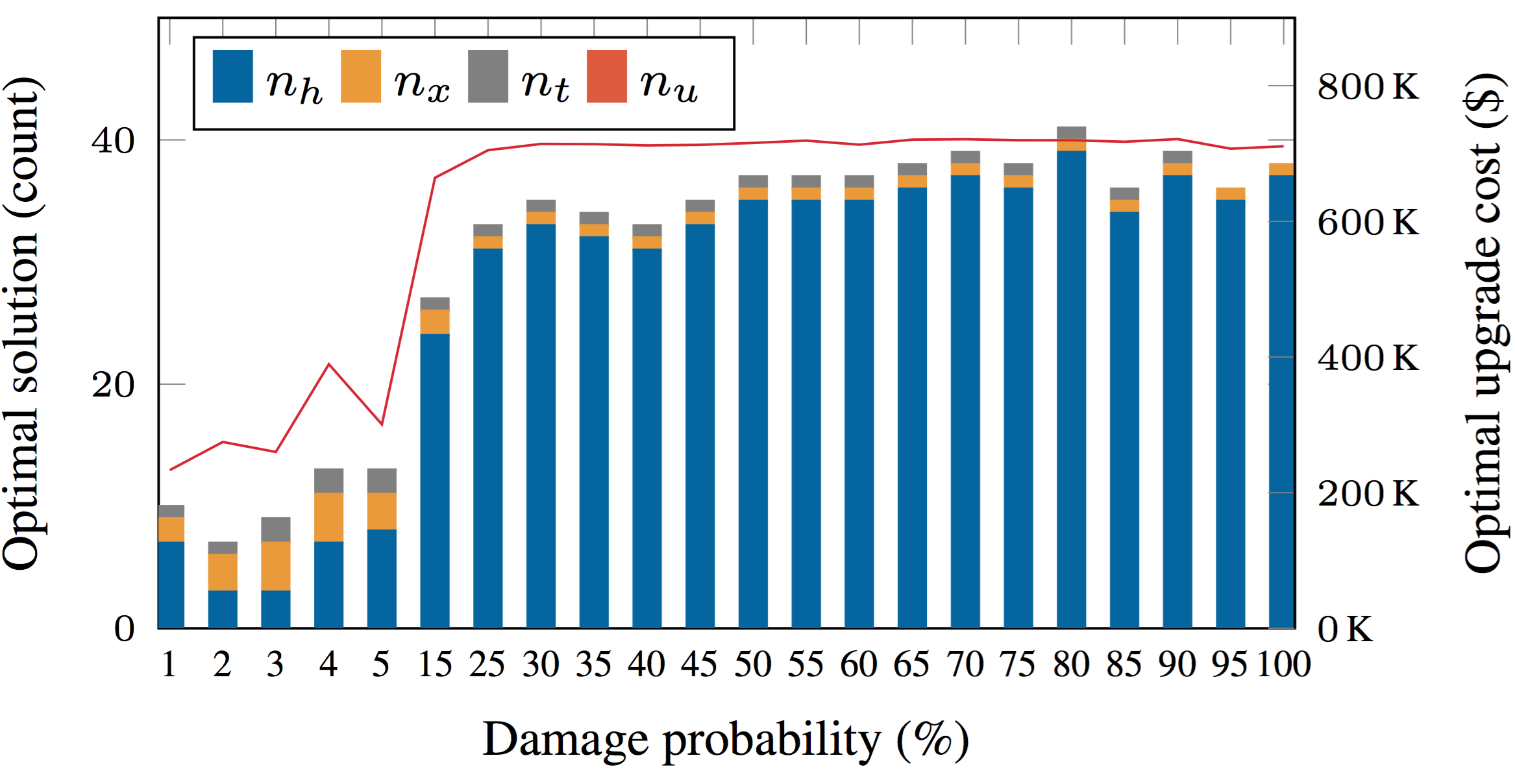}}\label{fig:opt_sol:urban}
%\end{figure}
%\begin{figure}[ht!] \ContinuedFloat
	\subfloat[Network123]{\includegraphics[width=0.5\textwidth]{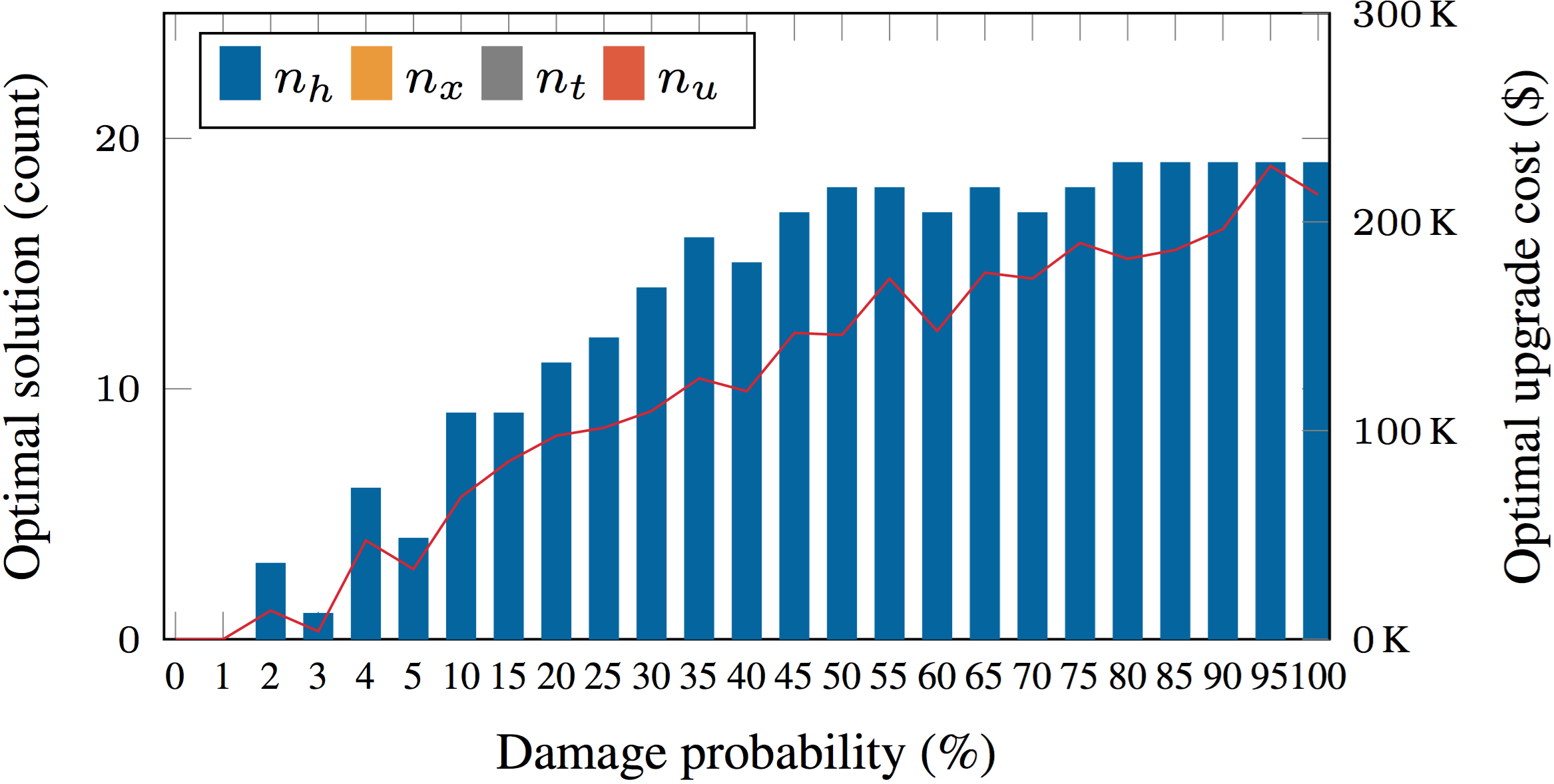}}\label{fig:opt_sol:network123}\caption{Statistics on the Optimal Grid Designs.} \label{fig:opt}
\end{figure}

Let $n_h, n_x, n_t,$ and $n_u$ be the number of hardened lines, new
lines, new switches, and new generators in the optimal design. Figure
\ref{fig:opt} reports these values for various damage levels and
the three networks. The red line indicates the optimal upgrade costs,
and the counts of the upgrade options are represented as a bar.  The
results show that hardening lines is the major component of each
optimal design and that its share increases with the disaster
intensity.  The results also show that DGs are used in significant
numbers in the rural network, while new lines and switches
complement hardening in the urban model. This was expected given the
length of the lines in these two networks. The third network only
needs line hardenings. 

\begin{figure}[!t]
\centering
\includegraphics[width=0.7\textwidth]{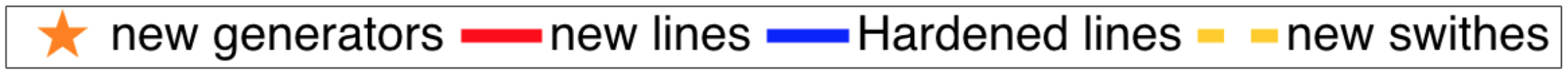}\\
\subfloat[Rural network] {\includegraphics[width=0.35\textwidth]{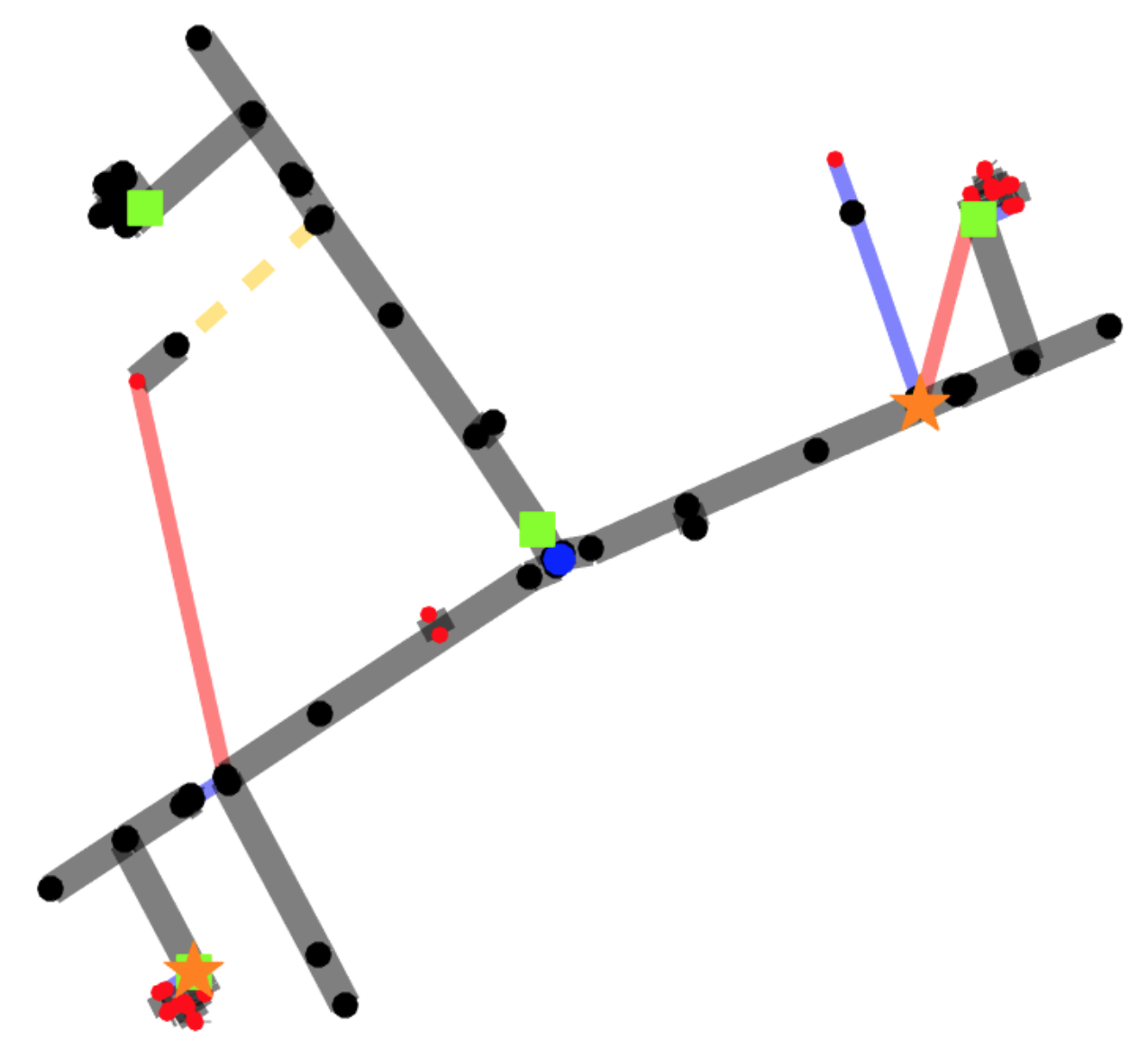}}\label{fig:opt_design:rural}
\subfloat[Urban network]
{\includegraphics[width=0.4\textwidth]{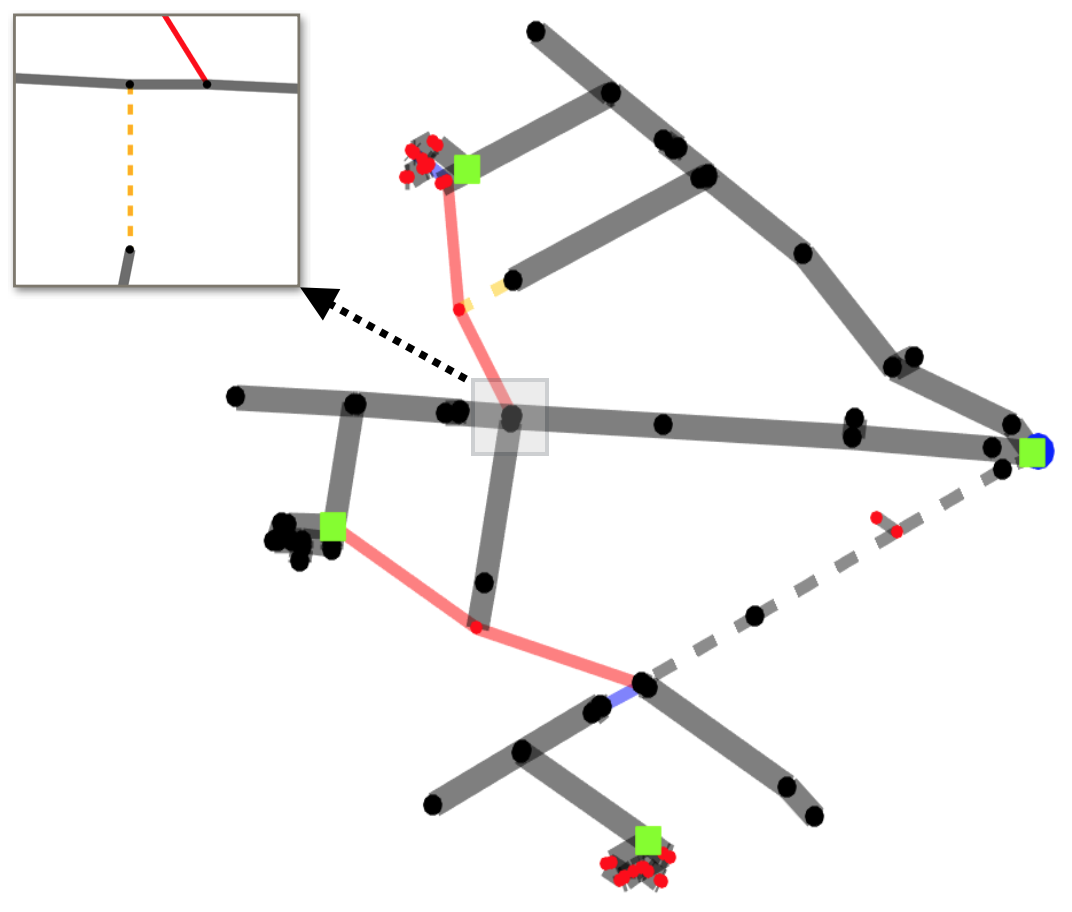}}\label{fig:opt_design:urban}
%\subfloat[Network123]
%{\includegraphics[width=0.45\textwidth]{Figures/optimal_design_network123.png}}\label{fig:opt_design:network123}
\caption{Optimal Designs of the Rural and Urban Networks (3\% damage level).}
\label{fig:opt_design}
\end{figure}

\begin{figure}[!t]
\centering
\includegraphics[width=0.6\textwidth]{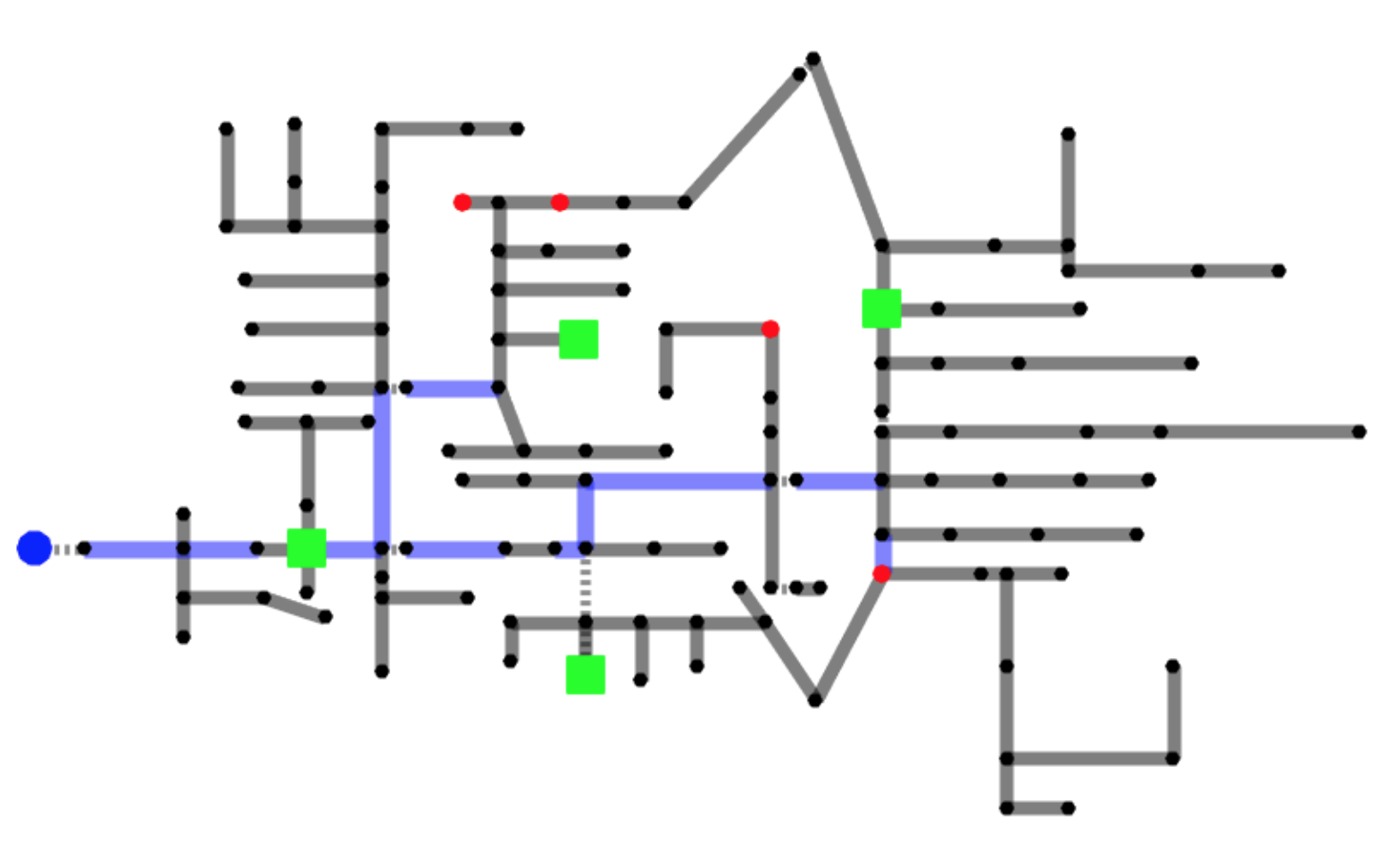}
\caption{Optimal Design of Network {\sc network123} (20\% damage level).}
\label{fig:opt_design:network123}
\end{figure}

\subsection{Impact of the Communication Network}
\label{sec:impact_commu}

First note that ignoring the communication network is equivalent to
assuming that every bus has its own communication center. In the
following, $\Gc(k)$ denotes a communication network with $k$ centers
and $\Gc(\infty)$ the case where each bus has a center.

\begin{figure}[!t]
\centering
\includegraphics[width = 0.35\textwidth]{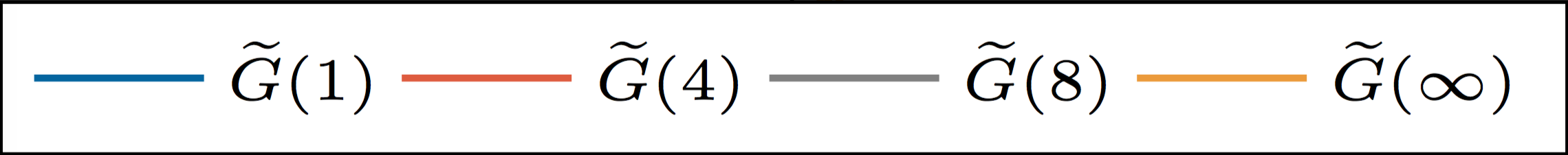}\\
\subfloat[Rural network] {\includegraphics[width=0.25\textwidth]{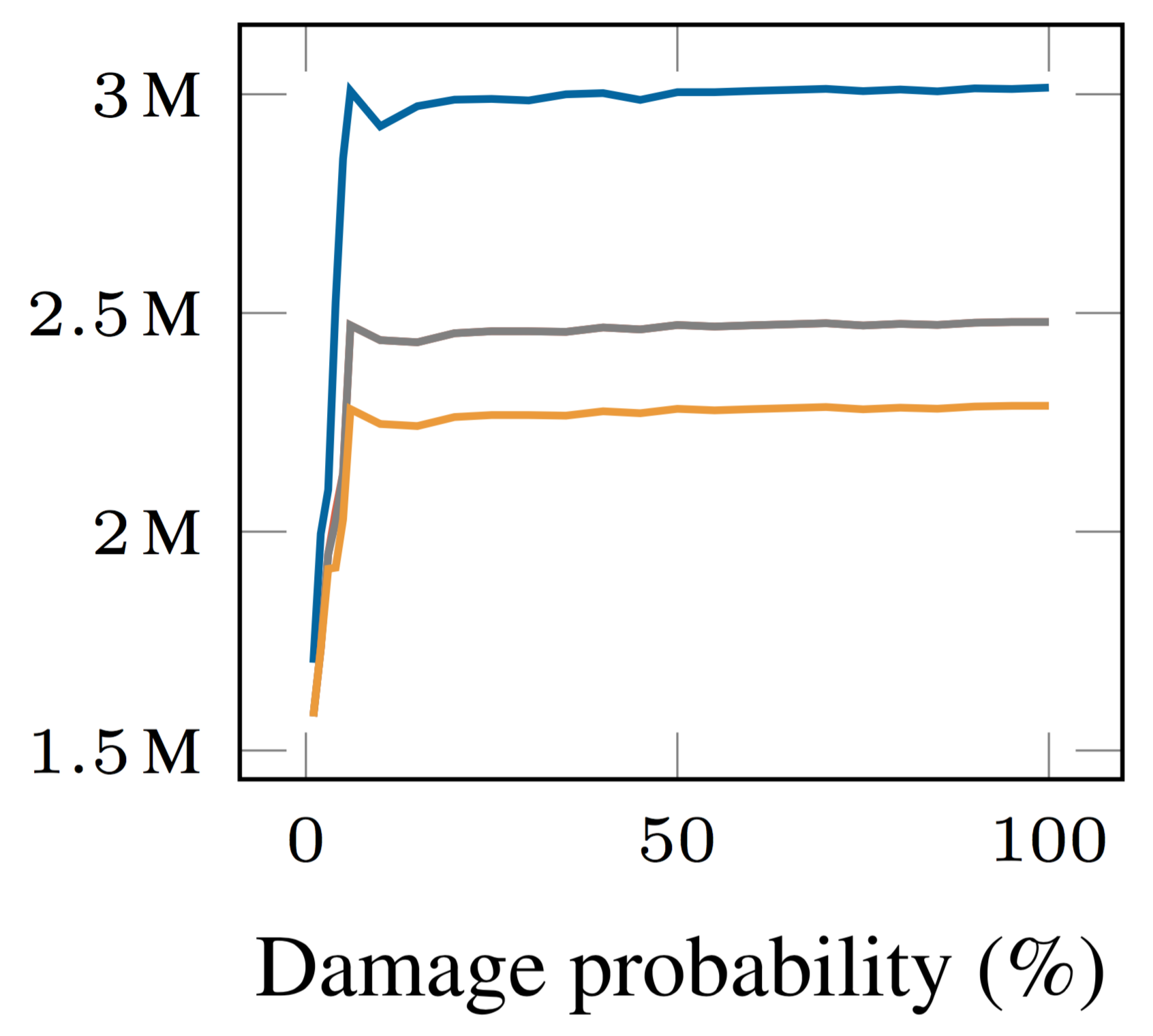}}\label{fig:cost_cc:rural}
\subfloat[Urban network]{\includegraphics[width=0.25\textwidth]{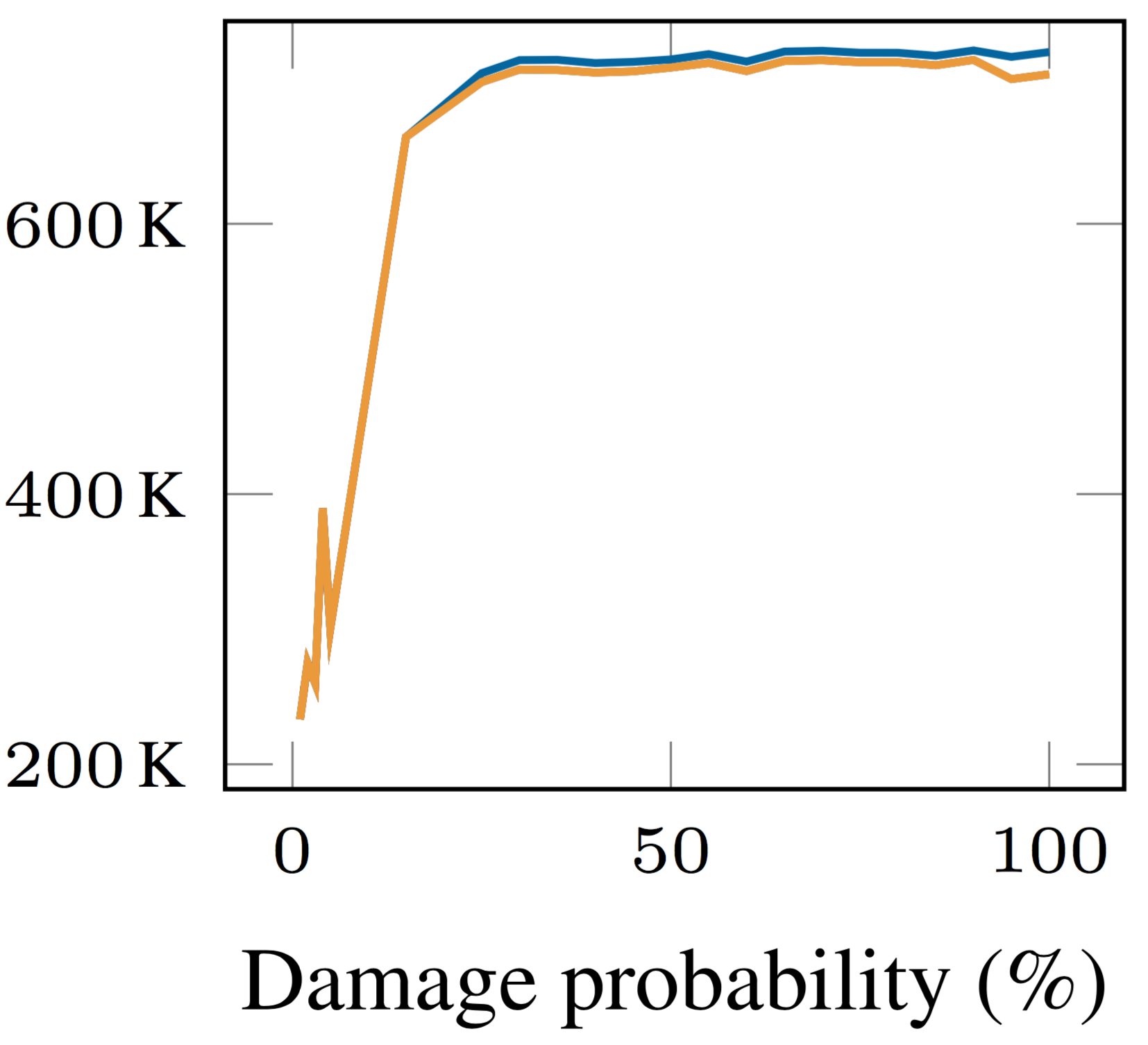}}\label{fig:cost_cc:urban}
\subfloat[Network123]{\includegraphics[width=0.25\textwidth]{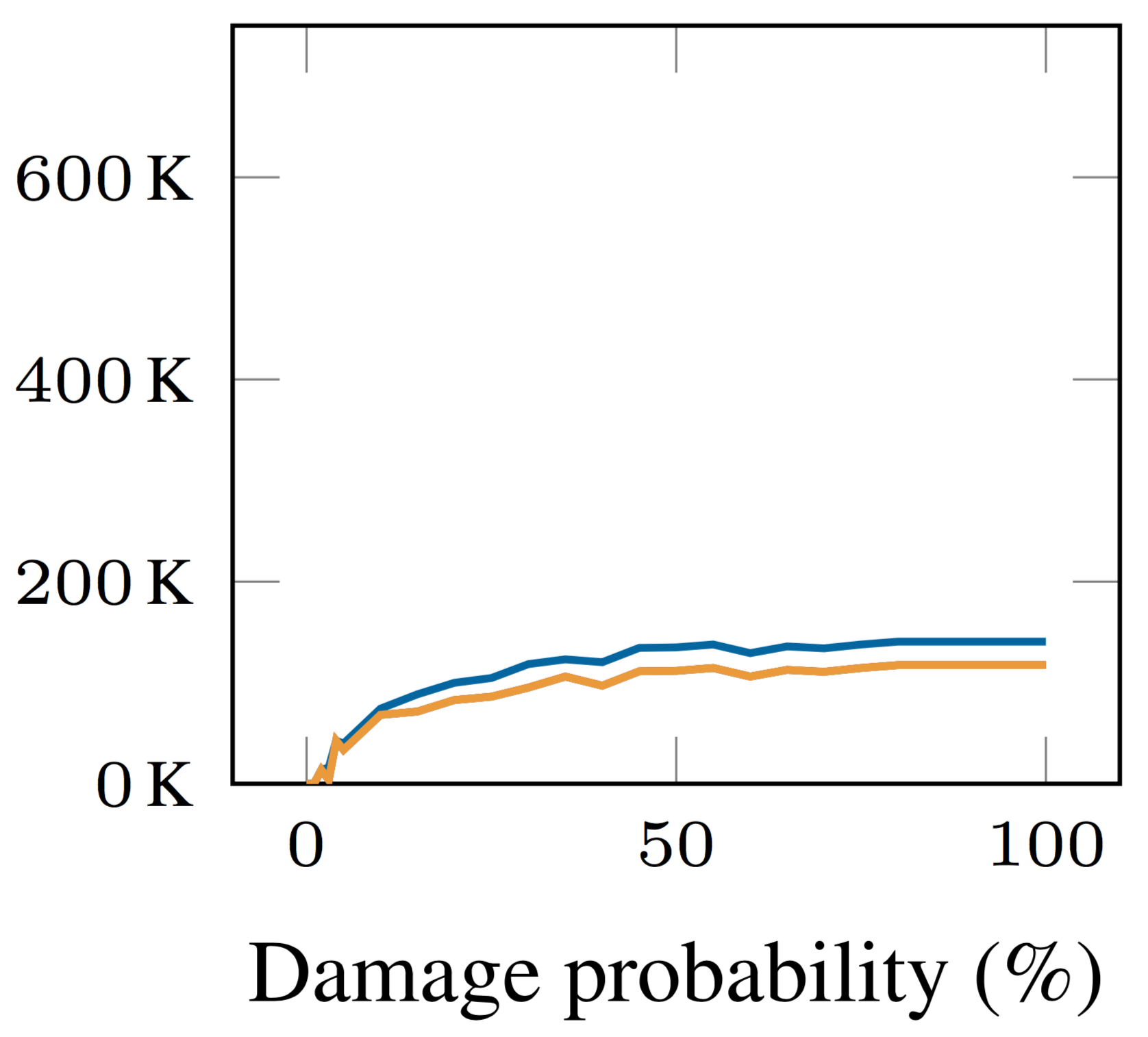}}\label{fig:cost_cc:network123}
\caption{Cost Analysis For the Number of Communication Centers}
\label{fig:cost_cc}
\end{figure}

Figure \ref{fig:cost_cc} and Table \ref{table:sol_cc} report the
impact of the communication system: They report optimal objective
values and solution statistics under various numbers of communication
centers. Fewer communication centers lead to significant
cost increases in the rural network, but have limited effect on 
the urban network and {\sc network123}. In the rural network,
resiliency comes from forming microgrids with DGs, which require
their own communication centers. When these are not available, 
optimal designs harden existing lines and build new lines and switches,
which are more costly as substantiated in Table \ref{table:sol_cc}.

\begin{figure}[!t]
	\centering
	\subfloat[damage scenarios]{\includegraphics[width=0.7\textwidth]{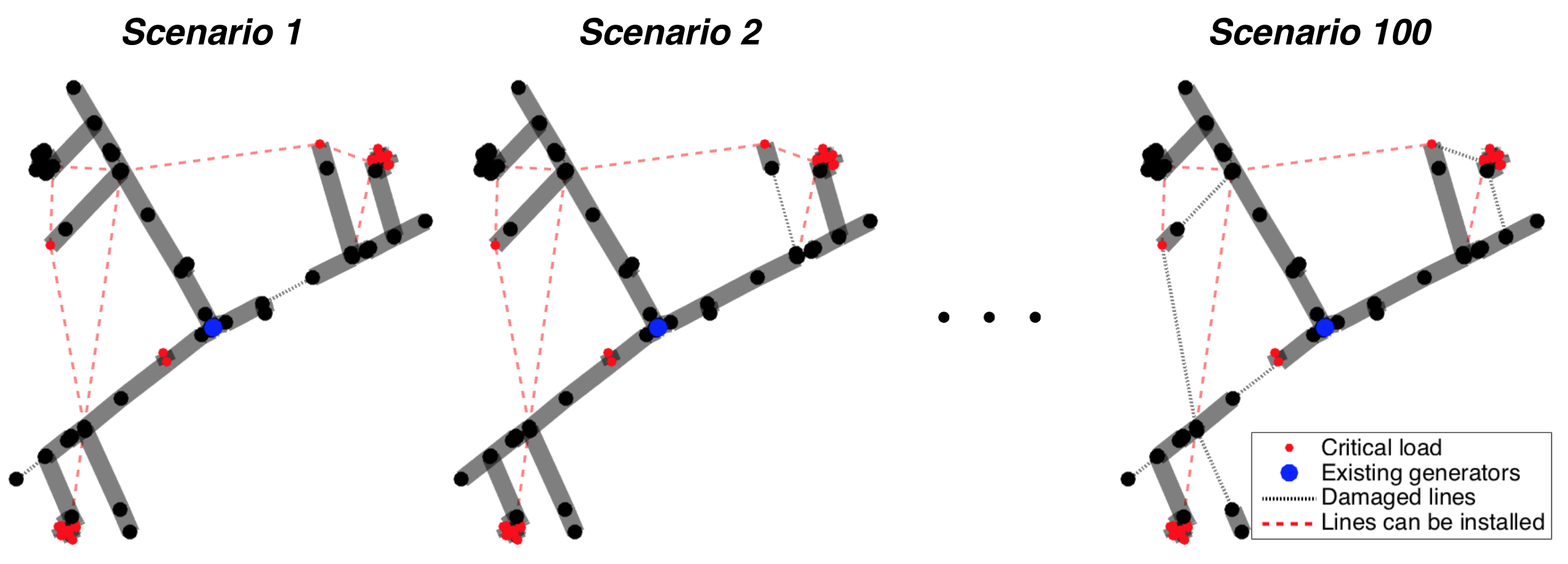}}\label{fig:opt_design:scenario}
	\subfloat[Optimal grid design]{\includegraphics[width=0.7\textwidth]{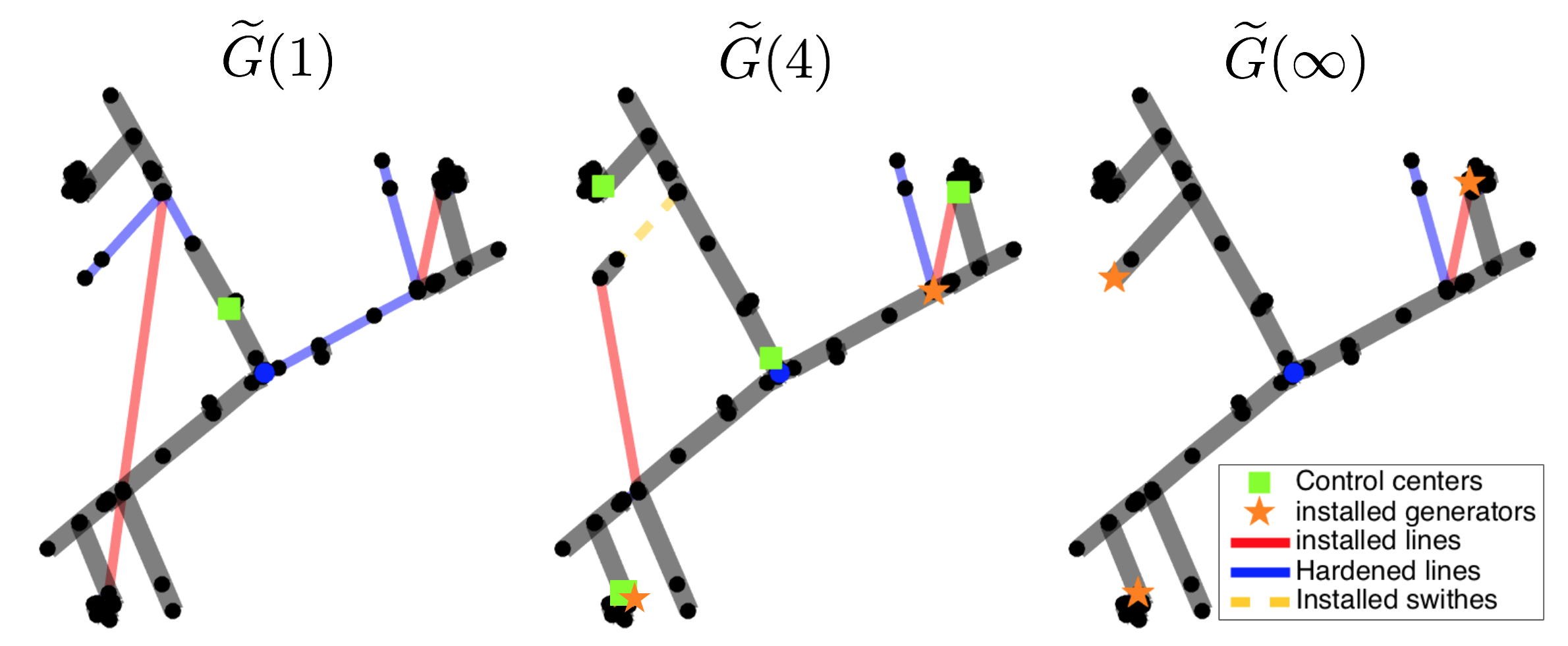}}\label{fig:opt_design:cc_1}
	\caption{Optimal Designs of the Rural Network under 3\% Damage and Various Communication Network Configurations.}
	\label{fig:optimal_design}
\end{figure}

Figure \ref{fig:optimal_design} illustrates the resulting designs on
the rural network for scenarios with a damage level of 3\%. The top
row depicts some of the scenarios and shows the affected lines. The
bottom row depicts the optimal designs for various configurations of
the communication network. For $\Gc(\infty)$, the optimal design
features three new DGs in the west-, north-, and east-end of the
network to meet the critical loads of each region. These regions are
then islanded under various scenarios. For $\Gc(4)$, the optimal
design installs a new line linking critical loads in the north side to
the west side of the network, instead of using DG in the north
side. This stems from Scenario 100 where a DG in the bus with critical
loads cannot be operated since it has no communication center.  For
$\Gc(1)$, scenario 1 prevents the operation of an east-end DG and
scenario 100 the operation of a west-end DG. Hence, the optimal design
only considers hardening and new lines and switches. On the other
hand, the urban network and {\sc Network123} achieve resiliency by
increasing grid connectivity for all communication networks.

\begin{table}[!t]
\TABLE{Impact of the Communication Network on Optimal Grid Designs. \label{table:sol_cc}}
{
%\caption{Impact of the Communication Network on Optimal Grid Designs}
%\label{table:sol_cc}
%\centering
\begin{tabular}{lcccccc}
\hline
%&  &				&	\multicolumn{4}{c}{number of} \\
%\cmidrule{4-7} 
& \up \down \makecell{Comm. Network} &	Obj. &	$n_h$	& $n_x$	& $n_t$	& $n_u$\\
\hline
%\hline
\up \multirow{4}{1cm}{Rural, 3\% damage}&	$\Gc(1)$&	2095.74 &	12	& 3	& 1	&0 \\
	&$\Gc(4)$&	1948.09	&6	&2&	1	&2\\
	&$\Gc(8)$&	1948.09	&6	&2	&1	&2\\
	\down &$\Gc(\infty)$&	1914.99&	5&	1&	0&	3\\
%	\midrule
%	\multirow{4}{1cm}{rural, 4\% damage}&	$\Gc(1)$	&2524.57&	26&	1&	1&	0\\
%	&$\Gc(4)$&	2046.80&	6&	2&	0&	3\\
%	&$\Gc(8)$&	2023.09&	5&	2&	0&	3\\
%	&$\Gc(\infty)$&	1917.7&	5&	1&	0&	3\\
%	\midrule
%\multirow{4}{1cm}{rural, 15\% damage} &	$\Gc(1)$	& 2972.69 &	34 &	1 &	1 &	0\\
%	&$\Gc(4)$& 	2433.05&	18&	1&	0&	3\\
%	&$\Gc(8)$&	2433.05&	18&	1&	0&	3\\
%	&$\Gc(\infty)$&	2241.45&	18&	0&	0&	3\\
	\hline
	\end{tabular}
    }
    {}
\end{table}

%The restrictive number of communication centers also affects the marginal cost to increase the resiliency level of the rural network. Figure \ref{fig:resiliency} illustrates the cost difference with regard to the resiliency level for each type of communication network. Note that $\Gc(1)$ not only has the highest total upgrade costs, but also the biggest marginal cost to enhance the resiliency level from 50\% to 80\%. Interestingly, when we increase the resiliency level, the upgrade cost for $\Gc(\infty)$ also increases and becomes the same as $\Gc(4)$ and $\Gc(8)$. That's because, with the high level of resiliency, $\eta_t$, it is inevitable to use line options to meet the non-critical load. This implies that with a high level of resiliency, there is an optimal number of communication centers that has the same upgrade cost as $\Gc(\infty)$; in this case study for the rural network, it is 4. 
%
%\begin{figure}[!ht]
%	\centering
%	{\includegraphics[width=0.35\textwidth]{Figures/legend_resiliency.png}}\\
%	\subfloat[$\Gc(1)$]{\includegraphics[width=0.24\textwidth]{Figures/resiliency_cc1.png}}\label{fig:resiliency:cc1}
%	\subfloat[$\Gc(4)$]{\includegraphics[width=0.24\textwidth]{Figures/resiliency_cc4.png}}\label{fig:resiliency:cc4}
%	\subfloat[$\Gc(8)$]{\includegraphics[width=0.24\textwidth]{Figures/resiliency_cc8.png}}\label{fig:resiliency:cc8}
%	\subfloat[$\Gc(\infty)$]{\includegraphics[width=0.24\textwidth]{Figures/resiliency_cc0.png}}\label{fig:resiliency:cc0}
%	\caption{Cost analysis on $\eta_t$, rural network }	\label{fig:resiliency}
%\end{figure}

\section{Performance Analysis of the Branch and Price Algorithm} 
\label{sec:performance}

This section studies the performance of the B\&P algorithm. All
computations were implemented with the C++/Gurobi 6.5.2 interface. They
use a Haswell architecture compute node configured with 24 cores (two
twelve-core 2.5 GHz Intel Xeon E5-2680v3 processors) and 128 GB RAM.
	
\subsection{ Computational Performance}
	
\begin{figure}[!t]
\centering
\subfloat[Computation Times]{\includegraphics[width = 0.4 \textwidth]{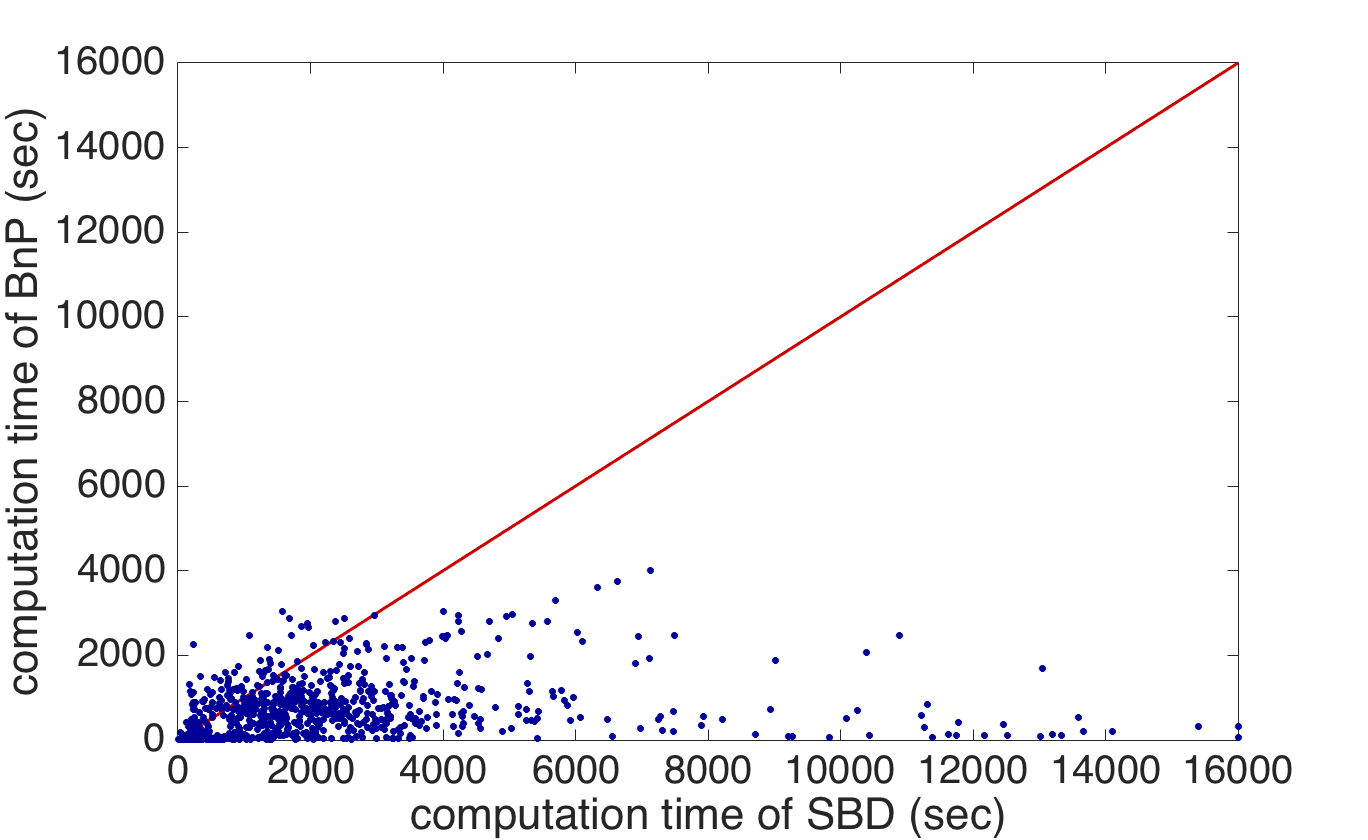}\label{fig:time_comparison:scatter}}
\subfloat[Average Computation Times and 95\% Confidence Intervals]{\includegraphics[width = 0.45 \textwidth]{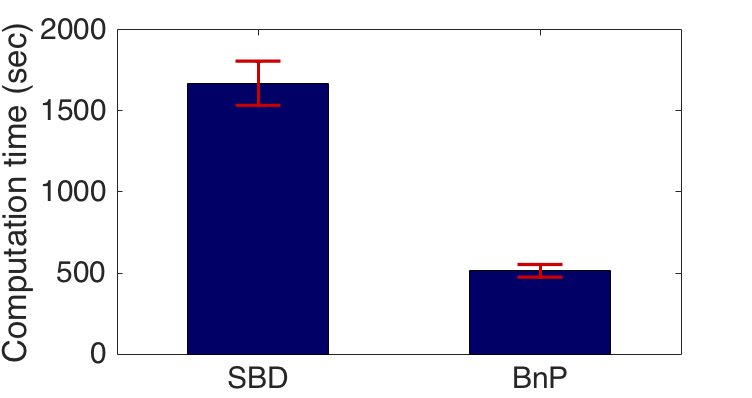}\label{fig:time_comparison:bar}}
\caption{Comparison of Computation Times: SBD versus B\&P.}	
\label{fig:time_comparison}
\end{figure}

Figure \ref{fig:time_comparison:scatter} reports the computation time
of the B\&P and SBD algorithms for all the instances described in
Section \ref{sec:data}, where the reference line (in red) serves to
delineate when an algorithm is faster than the other. Their statistics
are displayed in Figure \ref{fig:time_comparison:bar}. In average, the
B\&P algorithm is faster than the SBD algorithm by a factor of
3.25. These figures also indicate that the SBD algorithm has a high
degree of performance variance. This comes from the nature of the
scenario set $\mathcal S$. If $S$ contains a dominating scenario and
the scenario has low index in $\mathcal{S}$, then the SBD algorithm
solves the problem quickly. Otherwise, the SBD may need a large number
of iterations and the MIP model keeps growing in size with each
iteration. For 2 out of 1120 instances, the SBD algorithm times out
(wallclock time limit of 4 hours). On the other hand, the B\&P
algorithm is stable across all instances. The B\&P algorithm also has
the additional benefit that it produces improving feasible solutions
continuously. In contrast, the SBD algorithm only produces a feasible
solution at optimality. Finally, the B\&P algorithm appears more
stable numerically than the SBD algorithm. For 5 out of 1120
instances, the B\&P algorithm yields a better optimal solution than
the SBD algorithm as shown in Table \ref{table:opt_sol_gap}. Each such
solution was validated for feasibility.
	
\begin{table}[!t]
\TABLE{Numerical Stability of the B\&P Algorithm.\label{table:opt_sol_gap}}
{\begin{tabular}[h]{cccc}
\hline
\up \down \multirow{2}{2cm}{Instance}& \multicolumn{2}{c}{Opt. obj.val} & \\
\cline{2-3} 
\up   & SBD & B\&P	& Gap \\
\hline
\up \makecell{Rural, 30\% damage,
$\eta_t = 0.5,$ $\Gc(4)$} & 2458.49 & 2453.79 &	-0.19 \%\\
\makecell{Rural, 30\% damage, 
	$\eta_t = 0.6,$ $\Gc(4)$} & 2458.49 & 2453.79 &	-0.19 \%\\
\makecell{Rural, 30\% damage, 
	$\eta_t = 0.7,$ $\Gc(4)$} &2524.68 & 2519.98 &	-0.19 \%\\
\makecell{Rural, 30\% damage, 
	$\eta_t = 0.8,$ $\Gc(4)$} & 2572.31 & 2567.60 &	-0.19 \%\\
\down \makecell{Network123, 55\% damage, 
	$\eta_t = 0.8,$ $\Gc(8)$} & 232.48 & 227.27 &	-2.24 \%\\
			\hline
\end{tabular}}
{}
\end{table}

\subsection{Solution Quality at the Root Node.}
\label{sec:result:tree}

The problem reformulation produces a strong lower bound and the
majority of the instances are proven optimal at the root node. Table
\ref{table:tree} summarizes the average number of branching nodes and
the average optimality gap at the root node. 
			
\begin{table}[!t]
\TABLE
{Branching Tree Statistics.\label{table:tree}}
{\begin{tabular}{cc}
\hline
\up \down Avg. \# of branching nodes &	Avg. opt. gap at the root node \\
\hline
\up \down 1.8 &	0.19 \%\\
\hline
\end{tabular}}
{}
\end{table}
	
\subsection{Benefits of the Accelerating Schemes}
\label{sec:result:accelerating}

To highlight its design choices, the B\&P algorithm is compared to a
column generation with dual stabilization \citep{du1999stabilized}.
In addition, the benefit of each of the accelerating schemes is
investigated independently by running the B\&P algorithm without the
considered extension. We sample 90 instances by setting $\eta_{t} =
0.5$ and $0.8$, and the damage level to $5\%,$ $30 \%,$ $65\%,$
$85\%,$ $100\%$ for the three networks $\Gc(0)$, $\Gc(1),$ and
$\Gc(4)$. Dual stabilization prevents dual variables from fluctuating
too much, which is often the case in column generation. It tries to
confine dual variables in a box that contains the current best
estimate of the optimal dual solution and penalizes solutions that
deviate from the box. See, for instance, \citet{du1999stabilized,
  lubbecke2005selected} for details about stabilized column
generation. Our implementation updates the box whenever the Lagrangian
lower bound is updated.

Table \ref{table:dual} summarizes the computational performance of the
stabilized column generation in comparison with the B\&P
algorithm. B\&P$_B$ denotes the branch-and-price algorithm with the
basic scheme only (Section \ref{sec:cg:basic}) and B\&P$_S$ stands for
the branch and price algorithm with dual stabilization. The symbol
$\dag$ is used to denote that the algorithm reaches the wallclock time
limit for some instances. For more than one third of the sampled
instances, B\&P$_B$ and B\&P$_S$ exceed the wallclock time limit. For
instances where both algorithms terminate within the time limit,
B\&P$_S$ is faster than B\&P$_B$ by a factor of around 4. Although the
dual stabilization does improve the computation time of the basic
algorithm, it is still not adequate to solve the ORDPDC
practically. The B\&P algorithm, on the other hand, shortens
computation times by a factor of 26.35.
	
\begin{table}[!t]
\TABLE
{Comparison to a Column Generation with Dual Stabilization. 
\label{table:dual}}
{\begin{tabular}{lrr}
\hline
& \up \down Avg. computation time (sec)  &	Avg. number of  iterations \\
\hline
\up B\&P$_{B}$&	12857.97$^\dag$	&3122.57$^\dag$	\\
B\&P$_{S}$ &	11563.44$^\dag$	&1514.58$^\dag$\\
\down B\&P&488.03	& 96.12\\
\hline
\end{tabular}}{}
\end{table}

\begin{figure}[!t]
\centering
\subfloat[Optimality cut] {\includegraphics[width = 0.33 \textwidth]{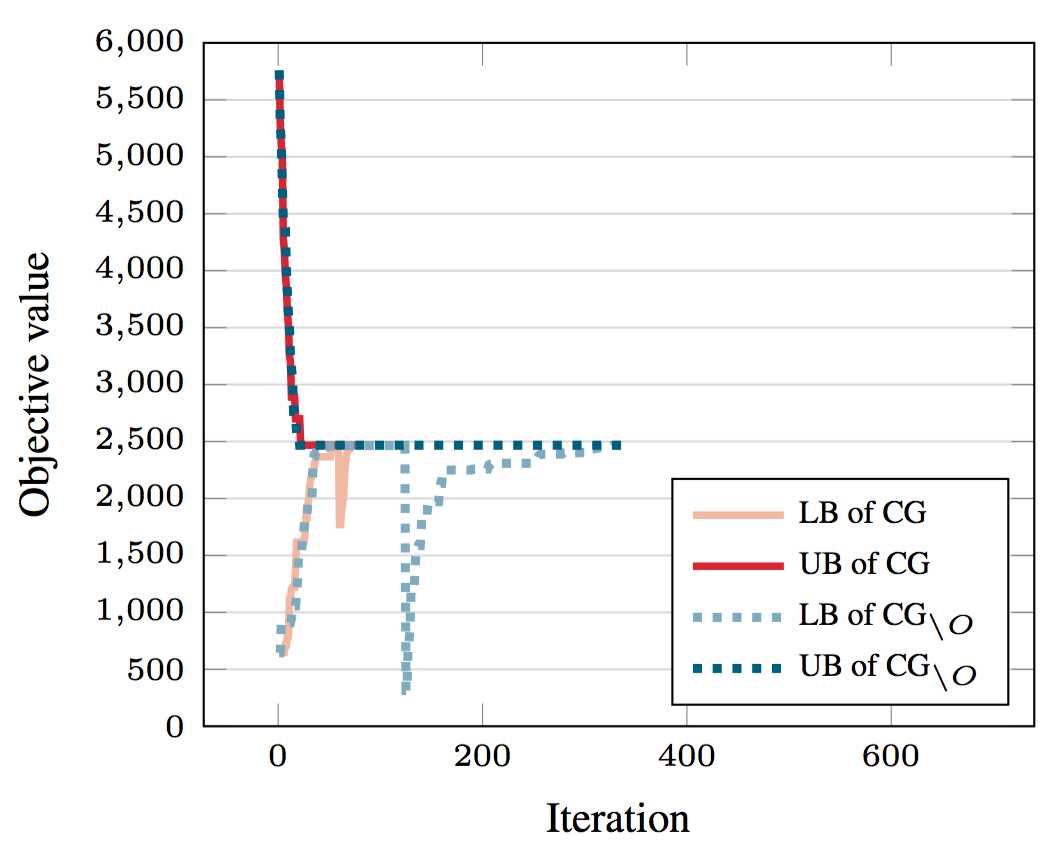}
\label{fig:convergence:O}}
\subfloat[	Revised reduced cost ] {
	\includegraphics[width = 0.33 \textwidth]{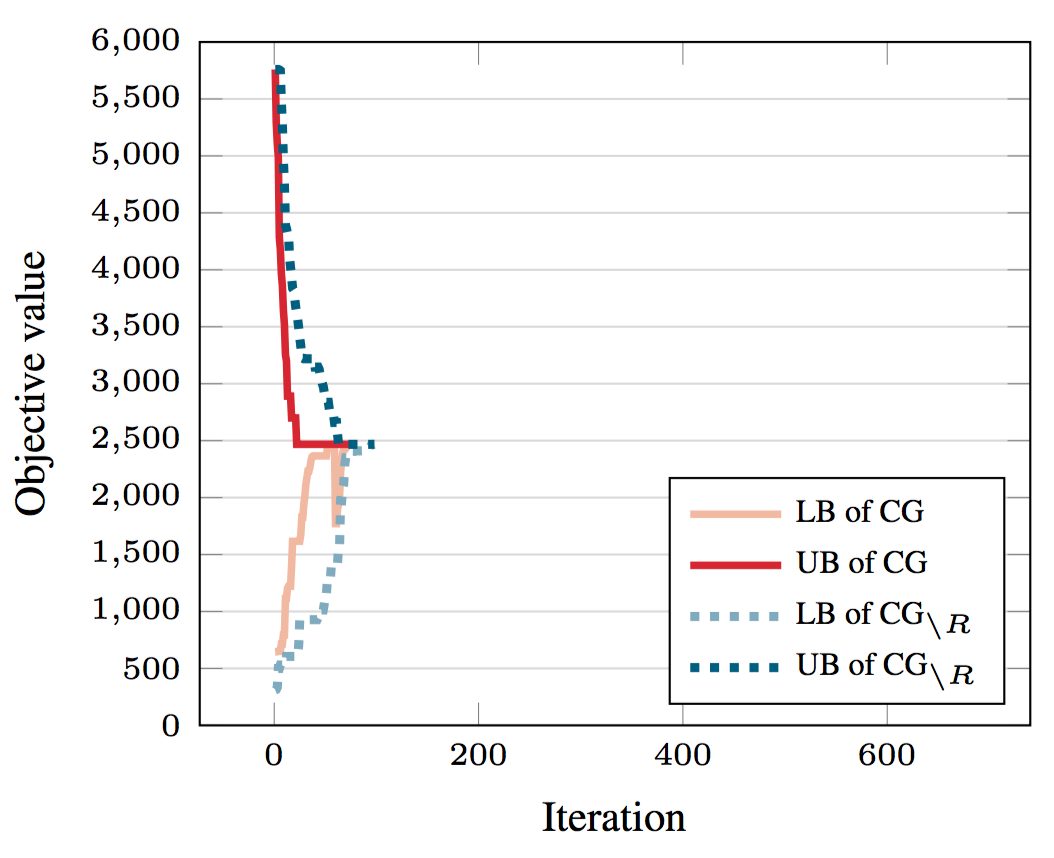}
\label{fig:convergence:R}}
	\subfloat[lexicographic objective function] {
		\includegraphics[width = 0.33 \textwidth]{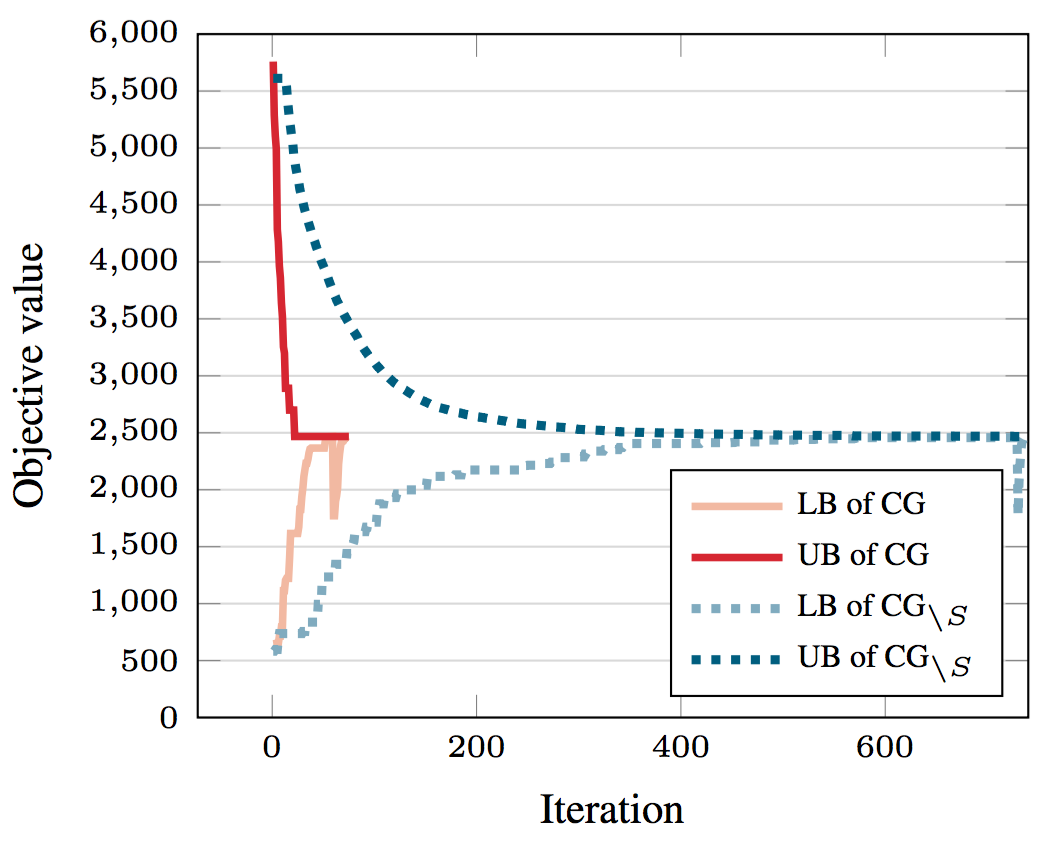}
\label{fig:convergence:S}}
\caption{Comparison of Convergence Rates (rural network, 6\% level of damage).}
\label{fig:convergence}
\end{figure}

The next results investigate the performance gain of each accelerating
scheme by removing them one at a time from the B\&P algorithm.  Table
\ref{table:CG:comparison} describes the computational performance and
Figure \ref{fig:convergence} illustrates the impact of each
accelerating schemes on the convergence rate of the rural network
under 6\% damage level. In the table and figure, $R$ denotes the
revised reduced cost, $C$ the optimality cut, $O$ the lexicographic objective
pricing problem, B\&P$_{\setminus k}$ the B\&P algorithm without
scheme $k$, with $k \in \{R, C, O\}$, and CG$_{\setminus k}$ the
column generation of B\&P without scheme $k$.

The results in Table \ref{table:CG:comparison} indicate that all the
accelerating schemes contribute to the computational performance of
the B\&P algorithm.  Figure \ref{fig:convergence:O} illustrates the
key role of the optimality cut. Without this cut, the second stage of
the column generation which uses the traditional pricing objective
does not take advantage of the columns generated in the first stage
and its lower bound drastically drops. Figure \ref{fig:convergence:R}
compares the convergence behavior of CG and CG$_{\setminus R}$,
showing that CG reaches the optimal objective value faster than
CG$_{\setminus R}$. Figure \ref{fig:convergence:S} highlights the
impact of the lexicographic objective function and shows that it
significantly contributes to the fast convergence of the algorithm.
	
\begin{table}[!t]
\TABLE{Benefits of the Accelerating Schemes. \label{table:CG:comparison}}
{\begin{tabular}[h]{lrr}
\hline
&  \up \down Avg. computation time (sec)  &	\up \down Avg. number of iterations\\
\hline
\up B\&P&488.03	& 96.12\\
B\&P$_{\setminus R}$ &	844.24 & 96.39 \\
B\&P$_{\setminus C}$ &	2589.55$^\dag$& 215.94$^\dag$\\
\down B\&P$_{\setminus O}$ &	2979.84$^\dag$& 544.65$^\dag$\\
\hline
\end{tabular}}
{}
\end{table}

%% file: Conclusions.tex
\section{Conclusions} 
\label{sec:conclusions}

This paper proposed an expansion planning model to improve the
resiliency of distribution systems facing natural disasters. The
planning model considers the hardening of existing lines and the
addition of new lines, switches, and distributed generators that would
allow a subpart of the system to operate as a microgrid. The expansion
model uses a 3-phase model of the distribution system. In addition, it
also considers damages to the communication system which may prevent
generators and switches to be controlled remotely. The input of the
expansion model contains a set of damage scenarios, each of which
specifying how the disaster affects the distribution system.

The paper proposed a branch and price algorithm for this model where
the pricing subproblem generates new expansions for each damage
scenario. The branch and price uses a number of acceleration schemes
to address significant degeneracy in the model. They include
a new pricing objective, an optimality cut, and a multi-objective
function to encourage sparsity in the generated expansions. The
resulting branch and price algorithm significantly improves the
performance of a scenario-based decomposition algorithm and a branch
and price with a stabilized column generation. The case studies show
that optimal solutions strongly depend on the grid topology and the
sophistication of the communication network. In particular, the
results highlight the importance of distributed generation for rural
networks, which necessitates a resilient communication system. 

The acceleration techniques presented in this paper are not limited to
the electricity distribution grid planning problem; They can be used
on problems with similar structure, i.e, two-stage stochastic problems
with feasibility recourse.

Future work will be devoted to applying and scaling these techniques
to instances with thousands of components.

\section*{Acknowledgements}

The work was funded by the DOE-OE Smart Grid R\&D Program in the
Office of Electricity in the US Department of Energy through the Grid
Modernization Laboratory Consortium project, \textit{LPNORM}.  It was
carried out under the auspices of the NNSA of the U.S. DOE at LANL
under Contract No. DE-AC52-06NA25396. It was also partly supported by
National Science Foundation grant NSF-1638331.